\documentclass[11pt]{article}

\usepackage{color}
\usepackage{dsfont}
\usepackage{amssymb}
\usepackage{amsmath}
\usepackage{amsthm}
\usepackage{enumerate}
\usepackage{hyperref}

\newtheorem{Theorem}{Theorem}[part]

\newtheorem{Proposition}{Proposition}[part]

\newtheorem{Lemma}{Lemma}[part]

\newtheorem{Remark}{Remark}[part]
\newtheorem{Example}{Example}[part]

\def \N{\mathbb{N}}
\def \R{\mathbb{R}}

\def \E{\mathbb{E}}
\def \F{\mathbb{F}}
\def \G{\mathbb{G}}

\def \P{\mathbb{P}}
\def \Q{\mathbb{Q}}

\def \Ac{{\cal A}}
\def \Uc{{\cal U}}
\def \Bc{{\cal B}}

\def \Ec{{\cal E}}
\def \Fc{{\cal F}}
\def \Gc{{\cal G}}

\def \Pc{{\cal P}}

\def \Sc{{\cal S}}
\def \Tc{{\cal T}}

\def \Yc{{\cal Y}}

\def \Zc{{\cal Z}}

\def \1{{\mathds 1}}

\def \ni{\noindent}

\def \eps{\varepsilon}

\def \ep{\hbox{ }\hfill$\Box$}

\def\reff#1{{\rm(\ref{#1})}}

\def\beqs{\begin{eqnarray*}}
\def\enqs{\end{eqnarray*}}
\def\beq{\begin{eqnarray}}
\def\enq{\end{eqnarray}}

\newcommand{\nc}{\newcommand}
\nc{\esssup}{\mathop{\mathrm{ess\;sup}}}

\begin{document}
\title{Mean-variance hedging on  uncertain time horizon\\ in a market with a jump 
}
\author{Idris Kharroubi\footnote{The research of the author benefited from the support of the French ANR research grant LIQUIRISK.}\\ \footnotesize{CEREMADE, CNRS UMR 7534}, \\\footnotesize{Universit\'e Paris Dauphine}\\ 
\footnotesize{ \texttt{kharroubi @ ceremade.dauphine.fr}} 
  \and Thomas Lim\footnote{The research of the author benefited from the support of the ``Chaire Risque de Cr\' edit'', F\' ed\' eration Bancaire Fran\c caise.}\\ \footnotesize{Laboratoire d'Analyse et Probabilit\'es,}\\ \footnotesize{Universit\'e d'Evry}\footnotesize{ and ENSIIE,}\\ \footnotesize{\texttt{lim @ ensiie.fr}}
   \and Armand Ngoupeyou\\ \footnotesize{Laboratoire de Probabilit\'es et Mod\`eles Al\'eatoires,}\\ \footnotesize{Universit\'e Paris 7}\\ \footnotesize{\texttt{armand.ngoupeyou @ univ-paris-diderot.fr}}
}

\date{}

\maketitle
\begin{abstract}
In  this work, we  study the problem of mean-variance hedging with a random horizon $T\wedge \tau$, where $T$ is a deterministic constant and $\tau$ is a jump time of the underlying asset price process. 
We first formulate this problem as a stochastic control problem and relate it to a system of BSDEs with a jump. We then provide a verification theorem which gives the optimal strategy for the mean-variance hedging using the solution of the previous system of BSDEs. Finally, we prove that this system of BSDEs admits a solution via a decomposition approach coming from filtration enlargement theory. 
 \end{abstract}

\vspace{1cm}



\setcounter{section}{0}

\section{Introduction}
\setcounter{equation}{0} \setcounter{Assumption}{0}
\setcounter{Theorem}{0} \setcounter{Proposition}{0}
\setcounter{Corollary}{0} \setcounter{Lemma}{0}
\setcounter{Definition}{0} \setcounter{Remark}{0}

In most financial markets, the assumption that the market is complete fails to be true.  In particular, investors cannot always hedge the financial products that they are interested in.   One possible approach to deal with this problem is mean-variance hedging. That is, for a given financial product with terminal value $H$ at a  fixed horizon time $T$ and an initial capital $x$, we need to find a strategy $\pi^*$ such that the value $V^{x,\pi^*}$ of the portfolio with initial amount $x$ and strategy $\pi^*$ minimizes the mean square error 
\beqs
\E\Big[\big|V^{x,\pi}_T-H\big|^2\Big]
\enqs
 over all possible investment strategies $\pi$. 

In this paper, we are concerned with the mean-variance hedging problem over a random horizon. More precisely, we consider a random time $\tau$ and a contingent claim with a gain at time $T\wedge\tau$ of the form 
\beq\label{form H intro}
H & = & H^b\mathds{1}_{T< \tau}+ H^a_\tau\mathds{1}_{T\geq \tau},
\enq
where $T <\infty$ is a fixed deterministic terminal time. We then study the mean-variance hedging problem over the horizon $[0,T\wedge\tau]$ defined by 
\beq\label{mevhdefintro}
\inf_{\pi}\E\Big[\big|V^{x,\pi}_{T\wedge \tau}-H\big|^2\Big]\;.
\enq
 Financial products with gains of the form \reff{form H intro} naturally appear on financial markets, see e.g.  Examples \ref{exle decom h assur}, \ref{exle decom h CDS} and \ref{exle decom h credit} presented in Subsection \ref{MVH pres}. 

The mean-variance hedging problem with deterministic horizon $T$ is  one of the classical problems from mathematical finance and has been considered by several authors via two main approaches.
One of them is based on martingale theory and projection arguments and the other considers the problem as a quadratic stochastic control problem and describes the solution using BSDE theory.
 
 The bulk of the literature primarily focuses on the continuous case where both approaches are used (see e.g. Delbaen and Schachermayer \cite{dealbenschach96},  Gouri\'eroux {\it et al.} \cite{GLP98}, Laurent and Pham \cite{LP99} and Schweizer  \cite{Schweiz96} for the first approach, and 
Lim and Zhou \cite{LZ02} and Lim \cite{L02} for the second one).

  In the discontinuous case, the mean-variance hedging problem is considered by Arai \cite{Ar05}, 
  Lim  \cite{limA09} and Jeanblanc {\it et al} \cite{JeMaSaSc11}.  In \cite{Ar05}, the author uses the projection approach for general semimartingale price processes model whereas in \cite{limA09} the problem is considered from the point of view of stochastic control for the case of diffusion price processes driven by Brownian motion and Poisson process. The author  provides under a so-called ``martingale condition" the existence of solutions to the associated BSDEs.  In the recent paper \cite{JeMaSaSc11}, the authors combine tools from both approaches, which allows them to work in a general semimartingale model and to give a description of the optimal solution to the mean-variance hedging via the BSDE theory.  More precisely the authors prove that the value process of the mean-variance hedging problem has a quadratic structure and that the coefficients appearing in this quadratic expression are related to some BSDEs. Then, they provide an equivalence between the existence of an optimal strategy and the existence of a solution to a BSDE associated to the control problem.  They have also shown in some specific examples, via the control problem,  the existence of solutions for BSDEs of interest. However the problem is still open in the general case.

In this paper, we study the mean-variance hedging with horizon $ T \wedge\tau$ given by \reff{mevhdefintro}. We use a stochastic control approach and describe the optimal solution by a solution to a system of BSDEs. 

We shall consider a model of diffusion price process driven by a Brownian motion and a random jump time $\tau$. We follow the progressive enlargement approach initiated by  Jacod, Jeulin and Yor  (see \cite{jeu80} and \cite{jeuyor85}), which leads to considering an enlargement of the initial information given by the Brownian motion to make $\tau$ a stopping time. 
We note that this approach allows to work under wide class of assumptions, in particular, on contrary to the Poisson case, no a priori law is fixed for the random time $\tau$. 

 Following the quadratic form obtained in \cite{JeMaSaSc11}, we use a  martingale optimality principle to obtain an associated system of nonstandard BSDEs. 
 We then establish a verification result (Theorem \ref{Th Verif}) which provides an explicit optimal investment strategy via the solution to the associated system of BSDEs. Our contribution is twofold.

\vspace{2mm}

\ni$\bullet$ We link the mean-variance hedging problem on a random horizon with a system of BSDE, in a general filtration progressive  enlargement setup which allows us to work without a priori knowledge of the law of jump part. We show that, under wide assumptions, the mean-variance hedging problem admits an optimal strategy described by the solution of the associated  system of BSDEs. 

\vspace{2mm}

\ni$\bullet$ We prove that the associated system of  BSDEs, which is nonstandard,  admits a solution. The main difficulty here is that the obtained system of BSDEs is nonstandard since it is driven by a  Brownian motion and a jump martingale and has generators with quadratic growth in the variable $z$ and are undefined for some values of the variable $y$. 
To solve these BSDEs we follow a decomposition approach inspired by the result of Jeulin (see Proposition \ref{propjeulin}) which allows to consider BSDEs in the smallest filtration (see Theorem \ref{Thm-exist-gen}). Then using BMO properties, we provide solutions to the decomposed BSDEs which lead to the existence of a solution to the BSDEs in the enlarged filtration.

\vspace{2mm}

We notice  that, for the problem at hand {\it i.e.} mean-variance hedging with horizon $T\wedge \tau$, the interest of our approach is that it provides a solution to the associated BSDEs,  without supposing any additional specific assumptions to the studied BSDEs unlike in \cite{limA09}  where to prove existence of a solution to the BSDE the author introduces the ``martingale condition''  or in \cite{JeMaSaSc11} where the existence of a solution to the BSDE is given in specific cases.

%

\vspace{2mm}

The paper is organized as follows. In Section 2, we present the details of the probabilistic model for the financial market, and setup the mean-variance hedging on random horizon.  In Section 3, we show how to construct the associated BSDEs via the martingale optimality principle and we state the two main theorems of this paper. The first one concerns the existence of a solution to the associated system of BSDEs and the second one is a verification theorem which gives an optimal strategy  via the solution of the BSDEs.  Then, Section 4 is dedicated to the proof of the existence of solution to the associated system of BSDEs. 
 Finally, some technical results are relegated to the appendix.
\section{Preliminaries and market model}
\label{modele}
\setcounter{equation}{0} \setcounter{Assumption}{0}
\setcounter{Theorem}{0} \setcounter{Proposition}{0}
\setcounter{Corollary}{0} \setcounter{Lemma}{0}
\setcounter{Definition}{0} \setcounter{Remark}{0}

\subsection{The probability space}

Let $(\Omega, \mathcal{G},\mathbb{P})$ be a complete probability space. 
We assume that this space is equipped with a one-dimensional standard Brownian motion $W$ and we denote by $\F:= (\Fc_t)_{t \geq 0}$ the right continuous complete filtration generated by $W$. We also consider on this space a random time $\tau$, which represents for example a default time in credit risk or in counterparty risk, or a death time in actuarial issues.  The random time $\tau$ is not assumed to be an $\F$-stopping time. We therefore  use in the sequel  the standard approach of filtration enlargement by considering $\G$ the smallest right continuous extension of $\F$ that turns $\tau$ into a $\G$-stopping time (see e.g. \cite{jeu80, jeuyor85}). More precisely $\G:=(\Gc_t)_{t\geq 0}$ is defined by
\beqs
\Gc_t & := & \bigcap_{\eps>0} \tilde \Gc_{t+\eps}\;,
\enqs 
for all $t\geq 0$, where $\tilde \Gc_{s} := \Fc_s\vee \sigma(\mathds{1}_{\tau\leq u}\;,u\in[0,s])$, for all $s\geq 0$.

\vspace{2mm}

We denote by $\Pc(\F)$ (resp. $\Pc(\G)$) the $\sigma$-algebra of $\F$ (resp. $\G$)-predictable subsets of $\Omega \times \R_{+}$, i.e. the $\sigma$-algebra generated by the left-continuous $\F$ (resp. $\G$)-adapted processes. 

\vspace{2mm}

We now introduce a decomposition result for $\Pc(\G)$-measurable processes. 
\begin{Proposition}\label{propjeulin}
Any $\Pc(\G)$-measurable process $X=(X_t)_{ t \geq 0}$ is represented as
\beqs
X_t & = & X^b_t \1_{t \leq \tau }+ X^a_t(\tau) \1_{t > \tau} \;,
\enqs
for all $ t \geq 0$, where $X^b$ is $\Pc(\F)$-measurable and $X^a$ is $\Pc(\F)\otimes\Bc(\R_+)$-measurable.
\end{Proposition}
{This result is proved in Lemma 4.4 of \cite{jeu80} for bounded processes and is easily extended to the case of unbounded processes. For the sake of completeness, we detail its proof in the appendix. 
}
\begin{Remark}\label{remjeulin}
{\rm  In the case where the studied process $X$ depends on another parameter $x$ evolving in a Borel subset $\mathcal{X}$ of $\R^p$, and if $X$ is $\Pc(\G)\otimes \Bc(\mathcal{X})$, 
  then, decomposition given by Proposition \ref{propjeulin} 
   is still true but where $X^b$ is $\Pc(\F)\otimes\Bc(\mathcal{X})$-mesurable and $X^a$ is $\Pc(\F)\otimes\Bc(\R_+)\otimes\Bc(\mathcal{X})$
   -measurable. Indeed, it is obvious for the processes generating $\Pc(\G)\otimes\Bc(\mathcal{X})$ 
   of the form $X_t(\omega, x)= L_t(\omega)R(x)$, $(t,\omega,x)$ $\in$ $\R_+\times\Omega\times\mathcal{X}$, where $L$ is $\Pc(\G)$
   -measurable and $R$ is $\Bc(\mathcal{X})$-measurable. Then, the result is extended to any $\Pc(\G)\otimes\Bc(\mathcal{X})$
   -measurable process by the monotone class theorem.
    }
\end{Remark}
We then impose the following assumption, which is classical in the filtration enlargement theory.

\vspace{2mm}

\noindent \textbf{(H)} The process $W$ remains a $\G$-Brownian motion. 

\vspace{2mm}

We notice that under \textbf{(H)}, the stochastic integral $\int_0^t X_sdW_s$ is well defined for all $\Pc(\G)$-measurable processes $X$ such that $\int_0^t|X_s|^2ds<\infty$. \\

\noindent In the sequel we denote by $N$ the process $\1_{\tau \leq .}$ and we suppose 

\vspace{2mm}

\noindent \textbf{(H$\tau$)} The process $N$ admits an $\F$-compensator  of the form $\int_0^{.\wedge \tau} \lambda_tdt$, i.e. $N-\int_0^{.\wedge \tau} \lambda_t dt$  is a $\G$-martingale, where $\lambda$ is a bounded $\Pc(\F)$-measurable process.

\vspace{2mm}
\noindent We then denote by $M$ the $\G$-martingale defined by
\beqs
M_t & := & N_t-\int_0^{t \wedge \tau}\lambda_sds\;,
\enqs
for all $t\geq 0$. We also introduce the process $\lambda^\G$ which is defined by $\lambda^\G_t:= (1-N_t)\lambda_t$.

\vspace{2mm}
\subsection{Financial model}
We consider a financial  market model  on the time interval $[0, T]$ where $0<T<\infty$ is a finite time horizon. We suppose that the financial market is composed by a riskless bond with zero interest rate and a risky asset  $S$. 
 The price process $(S_t)_{t \geq 0}$ of the risky asset is modeled by the linear stochastic differential equation
\beq\label{EDS S}
S_t & = & S_0+ \int_0^tS_{s^-} ( \mu_s ds + {\sigma}_s dW_s+ \beta_s dM_s)\;, \quad \forall t\in [0,T]\;,
\enq
where $\mu$, $\sigma$ and $\beta$  are $\Pc(\G)$-measurable processes and $S_0$ is a positive constant. For example $S$ could be a Credit Default Swap on the firm whose default time is $\tau$.
We impose the following assumptions on the coefficients $\mu$, $\sigma$ and $\beta$.

\vspace{2mm}

\ni \textbf{(H$S$)} \begin{enumerate}[(i)]
\item The processes $\mu$ and $\sigma$ are bounded: there exists a constant $C>0$ such that
\beqs
|\mu_t|+|\sigma_t| & \leq & C\;, \quad \forall t\in[0,T]\;,\quad \P-a.s. 
\enqs 
\item The process $\sigma$ is uniformly elliptic: there exists a constant $C>0$ such that
\beqs
|\sigma_t| & \geq & C\;, \quad \forall t\in[0,T]\;,\quad \P-a.s. 
\enqs 

\item There exists a constant $C$ such that
\beqs
-1~ \leq~ \beta_t~\leq~C\;,\quad \forall t\in[0,T]\;,\quad\P-a.s. 
\enqs
\end{enumerate}

\vspace{2mm}

\noindent Under \textbf{(H$S$)}, 
we know from e.g. Theorem 1 in \cite{Emery79} that the process $S$ defined by \reff{EDS S} is well defined.

\subsection{Mean-variance hedging}\label{MVH pres}
We consider investment strategies which are 
 $\Pc(\G)$-measurable processes $\pi$  such that 
\beqs
\int_0^{T\wedge\tau} |\pi_t|^2dt & < & +\infty\;,\qquad \P-a.s.
\enqs
This condition and \textbf{(H$S$)} ensure that the stochastic integral $\int_0^t{ \pi_r\over S_{r^-}} dS_r$ is well defined for such a strategy $\pi$ and $t\in [0,T \wedge \tau]$. 
The wealth process $V^{x,\pi}$ corresponding to a pair $(x,\pi)$, where $x\in \R$ is the initial amount, is defined by the stochastic integration 
\beqs
V^{x,\pi}_t & := & x+\int_0^t{\pi_r \over S_{r^-}}dS_r\;,\quad \forall \,t\in [0,T\wedge\tau]\;.
\enqs

We denote by $\Ac$ the set of admissible strategies $\pi$ such that
\beqs
\E\Big[\int_0^{T\wedge \tau} |\pi_t|^2dt\Big] & < & \infty\;.
\enqs
For $x \in \R$, 
 the problem of mean-variance hedging consists in computing the quantity
\beq\label{MVpb}
\inf_{\pi\in\Ac}\E \Big[\big|V_{T\wedge \tau}^{x,\pi}-H\big|^2 \Big]\;,
\enq
where $H$ is a bounded $\Gc_{T\wedge \tau}$-measurable random variable of the form
\beq\label{formeH}
H & = & H^b \mathds{1}_{T<\tau}+ H^a_\tau \mathds{1}_{T \geq \tau} \;,
\enq
where $H^b$ is an $\Fc_T$-measurable random variable valued in $\R$ and $H^a$ is a c\`ad-l\`ag $\Pc(\F)$-measurable process also valued in $\R$ and  such that
\beq\label{Habded}
\big\|H^b\big\|_{\infty}~<~\infty, & \mbox{ and } & \Big\|\sup_{t\in[0,T]} \big|H^a_t \big|\Big\|_{\infty}~<~\infty
\;,
\enq
where we recall that $\|.\|_\infty$ is defined by 
\beqs
{\|X\|}_\infty & := & \inf\Big\{ C\geq 0 ~:~\P\big(|X|\leq C\big)=1\Big\}\;,
\enqs 
for any random variable $X$. 
\vspace{2mm}

Since the problem we are interested in uses the values of the coefficients $\mu$, $\sigma$ and $\beta$ only on the interval $[0, T \wedge \tau]$, we can assume by Proposition \ref{propjeulin} that $\mu$, $\sigma$ and $\beta$ are $\Pc ( \F)$-measurable and we shall do that in the sequel.
\begin{Remark}
{\rm  For simplicity, we have supposed that the riskless interest rate is equal to zero. However, all the results can be extended to the case of a bounded $\Pc(\G)$-measurable interest rate process $r$. Indeed, for such an interest rate process the mean-variance hedging problem becomes
\beqs
\inf_{\pi\in\Ac}\E \Big[\big|\tilde V_{T\wedge \tau}^{x,\pi}-\tilde H\big|^2 \Big]\;,
\enqs
where $\tilde V^{x,\pi}$ and $\tilde H$ are the discounted values of $V^{x,\pi}$ and $H$ given by
\beqs
\tilde H & := & H \exp\Big(-\int_0^{T\wedge \tau}r_sds\Big)
\enqs
and 
\beqs
\tilde V^{x,\pi}_t & := & V^{x,\pi}_t \exp\Big(-\int_0^{t}r_sds\Big)\;,\quad t\in[0,T]\;.
\enqs
From the dynamic of $V^{x,\pi}$ we see that $\tilde V^{x,\pi}$ satisfies
\beqs
\tilde V^{x,\pi}_t & = & x+\int_0^t\pi_s\big(\tilde \mu_s ds+\tilde \sigma_sdW_s+\tilde \beta_sdM_s\big)
\enqs
where 
\beqs
\tilde \mu_t ~  := ~ e^{-\int_0^{t}r_sds}(\mu_t- r_t)\;,~\tilde \sigma_t~:=~e^{-\int_0^{t}r_sds}\sigma_t & \mbox{and} & \tilde \beta_t~:=~e^{-\int_0^{t}r_sds}\beta_t
\enqs
for $t\in[0,T]$. In particular, we get the same model but with coefficients $\tilde \mu$, $\tilde \sigma$ and $\tilde \beta$ instead of $\mu$, $\sigma$ and $\beta$. Since $\tilde \mu$, $\tilde \sigma$ and $\tilde \beta$ also satisfy (H$S$), we can extend the results to this model with new coefficients. 
}
\end{Remark}
\vspace{3mm}

\ni We end this section by two examples of financial products taking the form \reff{formeH}.

\begin{Example}[Insurance contract]\label{exle decom h assur} {\rm Consider a seller of an insurance policy which protects the buyer over the time horizon $[0,T]$ from some fixed loss $L$. Then if we denote by $\tau$ the time at which the loss appears, the losses of the seller are of the form
\beqs
H & = & -p \1_{T < \tau} +  (L-p)\mathds{1}_{T \geq \tau} \;,
\enqs
where $p$ denotes the premium that the insurance policy holder pays at time $0$.}
\end{Example}

\begin{Example}[Credit Default Swap with counterparty risk]\label{exle decom h CDS}
{\rm Consider a protection seller who sells a CDS against a credit event to a protection buyer for a nominal $N$ against a premium payments $p$ with a maturity $T$. If the reference entity defaults, the protection seller pays the buyer the nominal $N$ and the CDS contract is terminated. Moreover, both the buyer and seller of credit protection take on counterparty risk:
\begin{itemize}
\item the buyer takes the risk that the seller of credit protection may default, if the seller defaults the buyer loses its protection against default by the reference entity,
\item the seller takes the risk that the buyer may default on the contract, depriving the seller of the expected revenue stream.
\end{itemize}
Denote by $\tau$ the first default time, and by $\xi$ the random variable such that $\xi=1$ if the first default is the reference entity one and $\xi=0$  otherwise. The losses of the seller are of the form
\beqs
H &=& -pN T \1_{T < \tau} + N \1_{\tau \leq T, \xi =1} - pN \Big(\sum_{k=0}^{T}k\mathds{1}_{k\leq \tau<{k+1}}\Big) \1_{\tau \leq T}\;.
\enqs
}
\end{Example}
\begin{Example}[Credit contract]\label{exle decom h credit}
{\rm Consider a bank which lends an amount $A$ to a company over the period $[0,T]$. Suppose that the time horizon $[0,T]$ is divided on $n$ subintervals $[k{T\over n}, (k+1){T\over n}]$, $k=0,\ldots,n-1$, and that the interest rate of the loan over a time subinterval is $r$.  The company has then to pay ${(1+r)^n\over n}A$ to the bank at each time $k{T \over n}$, $k=1,\ldots,n$.  If we denote by $\tau$ the company default time, then the losses of the bank are given by 
\beqs
H & = & -((1+r)^n-1) A \mathds{1}_{T < \tau} + H^a_\tau\mathds{1}_{T\geq \tau} \;,
\enqs
where the function $H^a$ is given by 
\beqs
H^a _t & =  & -\sum_{k=1}^{n-1} \Big(k{(1+r)^n\over n }-1\Big)A\mathds{1}_{ k{T \over n} < t \leq (k+1){T \over n}}\;,\quad t \in[0,T]\,.
\enqs}
\end{Example}

\section{Solution of the mean-variance problem by BSDEs}
\subsection{Martingale optimality principle}\label{subsectionduale}\label{martingalprincopt}
To find the optimal value of the problem \reff{MVpb}, we follow the approach initiated by Hu \emph{et al.} \cite{huimkmul05} to solve the exponential utility maximization problem in the pure Brownian case. More precisely, we look for a family of processes 
\beqs
\Big\{\big( J^\pi_t \big)_{t\in[0,T]}~:~\pi\in\Ac\Big\}
\enqs
satisfying the following conditions
\begin{enumerate}[(i)]
\item $J^{\pi}_{T\wedge\tau}= {\big| V_{T\wedge\tau}^{x,\pi}-H \big|}^2$, for all $\pi\in \Ac$.

\item $J^{\pi_1}_0=J^{\pi_2}_0$, for all $\pi_1,\pi_2\in\Ac$. 

\item $\big( J^\pi_t \big)_{t\in[0,T]}$ is a $\G$-submartingale for all $\pi\in \Ac$.

\item  There exists some $\pi^*\in \Ac$ such that $\big( J^{\pi^*}_t \big)_{t\in[0,T]}$ is a $\G$-martingale.
\end{enumerate}

\noindent Under these conditions, we have 
\beqs
J_0^{\pi^*}  & = &   \inf_{\pi\in\Ac}\E \Big[\big|V_{T\wedge\tau}^{x,\pi}-H\big|^2 \Big]\;.
\enqs
Indeed, using (i), (iii) and Doob's optional stopping theorem, we have
\beq\label{inegHIM3}
J_0^\pi & \leq & \E\big[ J_{T \wedge \tau}^\pi \big]~=~\E\Big[ \big|V_{T\wedge\tau}^{x,\pi}-H\big|^2  \Big]\;,
\enq
for all $\pi\in \Ac$. Then, using (i), (iv)  and Doob's optional stopping theorem, we have 
\beq\label{inegHIM2}
J_0^{\pi^*} & = & \E\Big[ \big|V_{T\wedge\tau}^{x,\pi^*}-H\big|^2  \Big]\;.
\enq
Therefore, from (ii), \reff{inegHIM3} and \reff{inegHIM2}, we get for any $\pi \in \Ac$
\beqs
\E\Big[ \big|V_{T\wedge\tau}^{x,\pi^*}-H\big|^2  \Big] ~=~J_0^{\pi^*} & = & J_0^{\pi} ~\leq ~ \E\Big[ \big|V_{T\wedge\tau}^{x,\pi}-H\big|^2  \Big]\;.
\enqs
We can see that
\beqs
J_0^{\pi^*} & = & \inf_{\pi\in\Ac}\E\Big[ \big|V_{T\wedge\tau}^{x,\pi}-H\big|^2  \Big]\;.
\enqs

\subsection{Related BSDEs}
We now construct a family $\{(J_t^\pi)_{t\in[0,T]},~\pi\in \Ac\}$ satisfying the previous conditions by using BSDEs as in \cite{huimkmul05}.  To this end, we define the following spaces.
\begin{itemize}
\item $\Sc_{\G}^\infty$ is the subset of $\R$-valued c\`ad-l\`ag $\G$-adapted  processes $(Y_t)_{t\in[0,T]}$ essentially bounded
\beqs
{\| Y \|}_{\Sc^\infty} & := & \Big\|\sup_{t\in[0,T]}|Y_{t}|\Big\|_\infty~<~\infty \;.
\enqs
\item $\Sc_{\G}^{\infty,+}$  is the subset of $\Sc_{\G}^\infty$ of processes $(Y_t)_{t\in[0,T]}$ valued in $(0,\infty)$, such that
\beqs
{\Big\| {1\over Y} \Big\|}_{\Sc^\infty} & < & \infty \;.
\enqs
\item $L^2_{\G}$ is the subset of $\R$-valued $\Pc(\G)$-measurable processes $(Z_t)_{t\in[0,T]}$ such that
\beqs
\|Z\|_{L^2}  & := & \Big(\E\Big[ \int_0^T |Z_t|^2 dt \Big]\Big)^{1\over2} ~< ~\infty \;.
\enqs
\item $L^2(\lambda)$ is the subset of  $\R$-valued $\Pc(\G)$-measurable processes $(U_t)_{t\in [0,T]}$ such that
\beqs
\|U\|_{L^2(\lambda)} &:= & \Big(\E\Big[\int_0^{T \wedge \tau} \lambda_s|U_{s}|^2ds\Big]\Big)^{1\over2}~<~\infty\;.
\enqs
\end{itemize}

\noindent To construct a family $\{(J_t^\pi)_{t\in[0,T]},~\pi\in \Ac\}$ satisfying the previous conditions, we set 
\beqs
J^\pi_t & = & Y_t \big| V^{x,\pi}_{t\wedge\tau} - \Yc_t \big|^2 + \Upsilon_t\;,\quad t\in[0,T]\;,
\enqs
where\footnote{As commonly done for the integration w.r.t. jump processes, the integral $\int_a^b$ stands for $\int_{(a,b]}$.} $(Y, Z, U)$ is solution in $ \Sc^{\infty,+}_\G\times L^2_\G\times L^2(\lambda)$ to 
\beq\label{EDSRf^1}
Y_t & = & 1+ \int_{t\wedge \tau}^{T\wedge \tau}\mathfrak{f}(s, Y_s, Z_s, U_s)ds -\int_{t\wedge \tau}^{T\wedge \tau} Z_s dW_s -\int_{t\wedge \tau}^{T\wedge \tau} U_s dM_s\;,\quad t\in[0,T]\;,\qquad
\enq
 $(\Yc,\Zc,\Uc)$  is solution in $\Sc^\infty_\G\times L^2_\G\times L^2(\lambda)$ to
\beq\label{EDSRf^2}
\Yc_t & = & H+ \int_{t\wedge \tau}^{T\wedge \tau}\mathfrak{g}(s, \Yc_s, \Zc_s, \Uc_s)ds -\int_{t\wedge \tau}^{T\wedge \tau} \Zc_s dW_s -\int_{t\wedge \tau}^{T\wedge \tau} \Uc_s dM_s\;,\quad t\in[0,T]\;, \qquad
\enq
and $(\Upsilon,\Xi,\Theta)$  is solution in $\Sc^\infty_\G\times L^2_\G\times L^2(\lambda)$ to
\beq\label{EDSRf^3}
 \Upsilon_t & = & \int_{t\wedge \tau}^{T\wedge \tau}\mathfrak{h}(s, \Upsilon_s, \Xi_s, \Theta_s)ds - \int_{t\wedge \tau}^{T\wedge \tau}\Xi_s dW_s - \int_{t\wedge \tau}^{T\wedge \tau}\Theta_s dM_s \;,\quad t\in[0,T]\;. \qquad
\enq
\begin{Remark}\label{remUbdd}
{\rm We notice that the jump components $U$, $\Uc$ and $\Theta$ are also bounded since $Y$, $\Yc$ and $\Upsilon$ are in $\Sc^\infty_\G$. Indeed, let $C$ be a constant such that 
\beq\label{borne Y}
{\| Y \|}_{\Sc^\infty} & \leq & C \;.
\enq 
Then since $Y_{.-}+U$ is $\G$-predictable, we have
\beqs
\E\Big[\int_0^{T}\mathds{1}_{|Y_{t^-}+U_t|>C}\lambda_t^\G dt\Big] & = & \E\Big[\int_0^{T}\mathds{1}_{|Y_{t^-}+U_t|>C}d N_t\Big] \\
 & = & \E\Big[ \mathds{1}_{|Y_{\tau^-}+U_\tau|>C, \tau\leq T} \Big]\\
  & = & \E\Big[ \mathds{1}_{|Y_{\tau}|>C, \tau\leq T} \Big]\\
   & = & 0\;.
\enqs
Therefore, we have $|Y_{.^-}+U|\leq C$  in $L^2(\lambda)$. From \reff{borne Y} we get $|U|\leq 2C$ in $L^2(\lambda)$.  The same argument can be applied for $\Uc$ and $\Theta$.}  
\end{Remark}
In these terms, we are bound to choose three functions $\mathfrak{f}$, $\mathfrak{g}$ and $\mathfrak{h}$ for which $J^\pi$ is a submartingale for all $\pi \in \Ac$, and there exists a $\pi^* \in \Ac$ such that $J^{\pi^*}$ is a martingale. In order to calculate $\mathfrak{f}$, $\mathfrak{g}$ and $\mathfrak{h}$, we write $J^\pi$ as the sum of a (local) martingale $M^{\pi}$ and an (not strictly) increasing process $K^\pi$ that is constant for some $\pi^* \in \Ac$.\\
To alleviate the notation we write $\mathfrak{f}(t)$ (resp. $\mathfrak{g}(t)$, $\mathfrak{h}(t)$) for  $\mathfrak{f}(t,Y_t,Z_t,U_t)$ (resp. $\mathfrak{g}(t,\Yc_t,\Zc_t,\Uc_t)$, $\mathfrak{h}(t,\Upsilon_t,\Xi_t,\Theta_t)$) for $t\in [0,T]$.\\
Define  for each $\pi \in \Ac$ the process $X^\pi$ by 
\beqs
X^\pi_t & := & V^{x,\pi}_{t\wedge \tau} - \Yc_t\;,\quad t\in [0,T]\;.
\enqs 
From It\^o's formula, we get 
\beq\label{decompJpi}
dJ^\pi_t & = & dM^\pi_t + d K^\pi_t  \;, 
\enq
where $M^\pi$ and $K^\pi$ are defined by
\beqs
dM^\pi_t & := &  \Big\{ 2 X^\pi_{t^-} ( \pi_t \beta_t - \Uc_t ) ( Y_{t^-} + U_t) + | \pi_t \beta_t - \Uc_t |^2 ( Y_{t^-} + U_t) +{|X^\pi_{t^-}|^2U_t}+ \Theta_t \Big\} dM_t \nonumber \\
& & + \;\Big\{2Y_t X^\pi_t ( \pi_t \sigma_t - \Zc_t ) + Z_t |X^\pi_t|^2 + \Xi_t \Big\} dW_t \;, \label{def-M-pi}
\enqs
\beqs
dK^\pi_t & := & \Big\{Y_t \big[2 X^\pi_t ( \pi_t \mu_t +\mathfrak{g}(t)) + | \pi_t \sigma_t - \Zc_t |^2 \big] - |X^\pi_t|^2 \mathfrak{f}(t) + 2 X^\pi_t Z_t (\pi_t \sigma_t - \Zc_t)  \nonumber \\ && + \; 2 \lambda^\G_t X^\pi_t U_t (\pi_t \beta_t - \Uc_t)  + \lambda^\G_t  |\pi_t \beta_t - \Uc_t|^2 ( U_t + Y_t) - \mathfrak{h}(t) \Big\} dt \;.
\enqs
 We then write $dK^\pi$ in the following form
\beqs
dK^\pi_t & = & K_t(\pi_t)dt\;,
\enqs
where $K$ is defined by  
\beqs
K_t(\pi) & := & A_t |\pi|^2 + B_t \pi  + C_t\;, \quad \pi\in \R\;, ~\quad t\in[0,T]\;, 
\enqs
with 
\beqs
A_t & := & |\sigma_t|^2 Y_t + \lambda^\G_t|\beta_t|^2(U_t+Y_t)\;, \\
B_t & := & 2X^\pi_t(\mu_t Y_t +\sigma_tZ_t+\lambda^\G_t\beta_tU_t)-2\sigma_tY_t\Zc_t-2\lambda^\G_t\beta_t\Uc_t(Y_t+U_t)\;,\\
C_t & := & -\mathfrak{f}(t)|X^\pi_t|^2+ 2 X^\pi_t ( Y_t \mathfrak{g}(t) - Z_t \Zc_t-\lambda^\G_tU_t\Uc_t)+Y_t|\Zc_t|^2+\lambda^\G_t|\Uc_t|^2(U_t+Y_t)-\mathfrak{h}(t)\;,
\enqs
for all $t\in [0,T]$. To ensure that $K^\pi$ is nondecreasing  for any $\pi \in \Ac$ and that  $K^{\pi^*}$ is constant for some $\pi^*\in \Ac$, we take $K_t$ such that  $\min_{\pi \in \R} K_t(\pi) = 0$. Using $Y\in  \Sc^{\infty,+}_\G$ and \textbf{(H$S$)} (ii), we then notice that $A_t> 0$ for all $t\in[0,T]$. Indeed, we have
\beqs
0 & = & \E[[Y_{\tau}]^- \1_{\tau \leq T}] ~ = ~  \E[[Y_{\tau^-} +U_{\tau}]^- \1_{\tau \leq T}] ~ = ~
\E\Big[\int_0^T [Y_{s^-}+U_s]^-dN_s\Big] \;,
\enqs
therefore we get that 
\beqs
 \E\Big[\int_0^T [Y_{s^-}+U_s]^- dM_s + \int_0^T [Y_{s}+U_s]^-\lambda_s^\G ds\Big] &=& 0\;.
\enqs
From Remark \ref{remUbdd}, the  predictable process  $[Y_{.^-}+U]^-$ is bounded. Thus we get that the first integral on the left is a true martingale thus we have 
\beq\label{Y-+Upos}
  \E\Big[\int_0^T [Y_{s}+U_s]^-\lambda_s^\G ds\Big] &=& 0\;,
\enq
which gives $(Y_{s}+U_s)\lambda_s^\G\geq 0$ for $s\in[0,T]$. Therefore, the minimum of $K_t$ over $\pi\in \R$ is given by
\beqs
\underline K_t ~ := ~ \min_{\pi\in \R} K_t(\pi) & = & C_t-{|B_t|^2\over 4A_t}\;.
\enqs
We then obtain from the expressions of $A$, $B$ and $C$ that 
\beqs
\underline K _t & = & \mathfrak{A}_t |X_t^\pi|^2 + \mathfrak{B}_t X_t^\pi+ \mathfrak{C}_t\;,
\enqs
with
\beqs
\mathfrak{A}_t & := & -\mathfrak{f}(t)-\frac{|\mu_t Y_t+\sigma_t Z_t+\lambda^\G_t\beta_t U_t|^2}{|\sigma_t |^2Y_t + \lambda^\G_t |\beta_t|^2 (U_t + Y_t)}\;,\\
\mathfrak{B}_t & := & 2\Big\{\frac{(\mu_t Y_t+\sigma_t Z_t+\lambda^\G_t\beta_t U_t)(\lambda^\G_t\beta_t \Uc_t (Y_t+U_t)+\sigma_tY_t\Zc_t)}{|\sigma_t |^2Y_t + \lambda^\G_t |\beta_t|^2 (U_t + Y_t)}+\mathfrak{g}(t)Y_t-Z_t\Zc_t-\lambda^\G_tU_t\Uc_t\Big\}\;,\\
\mathfrak{C}_t & := & -\mathfrak{h}(t) + |\Zc_t|^2 Y_t + \lambda^\G_t (U_t + Y_t) |\Uc_t|^2
{-\frac{|\sigma_tY_t\Zc_t+\lambda^\G_t\beta_t\Uc_t(U_t+Y_t)|^2}{|\sigma_t|^2 Y_t + \lambda^\G_t |\beta_t|^2 ( U_t + Y_t)}}\;.\\
\enqs
For that the family $(J^\pi)_{\pi \in \Ac}$ satisfies the conditions (iii) and (iv) we choose $\mathfrak{f}$, $\mathfrak{g}$ and $\mathfrak{h}$ such that 
\beqs
\mathfrak{A}_t~=~0\;, ~\mathfrak{B}_t~=~0 & \text{and} & \mathfrak{C}_t~=~0\;,
\enqs
 for all $t\in [0,T]$. 
This leads to the following choice for the drivers $\mathfrak{f}$, $\mathfrak{g}$ and $\mathfrak{h}$
\begin{equation*}
\left\{\begin{array}{rcl}
\mathfrak{f}(t,y,z,u) & := & -\displaystyle{\frac{|\mu_t y + \sigma_t z + \lambda^\G_t \beta_t u|^2}{|\sigma_t |^2y + \lambda^\G_t |\beta_t|^2 (u + y)}}\;,\\ \\
\mathfrak{g}(t,y,z,u) & := & \frac{1}{Y_t} \Big[ Z_t z + \lambda^\G_t U_t u - \displaystyle{ \frac{(\mu_t Y_t + \sigma_t Z_t + \lambda^\G_t \beta_t U_t)(\sigma_t Y_t z + \lambda^\G_t \beta_t  (U_t + Y_t) u)}{|\sigma_t|^2 Y_t + \lambda^\G_t |\beta_t|^2 ( U_t + Y_t)} \Big]}\;,\\\\
\mathfrak{h}(t,y,z,u) & := & |\Zc_t|^2 Y_t + \lambda^\G_t (U_t + Y_t) |\Uc_t|^2
\displaystyle{-\frac{|\sigma_tY_t\Zc_t+\lambda^\G_t\beta_t\Uc_t(U_t+Y_t)|^2}{|\sigma_t|^2 Y_t + \lambda^\G_t |\beta_t|^2 ( U_t + Y_t)}}\;.
\end{array}\right.
\end{equation*}
We then notice that the obtained system of BSDEs is not fully coupled, which allows to study each BSDE alone as soon as we start from the BSDE $(\mathfrak{f},1)$\footnote{The notation BSDE $(f,H)$ holds for the BSDE with generator  $f$ and terminal condition  $H$.} and end with the BSDE $(\mathfrak{h},0)$.  However the obtained generators are nonstandard since they involve the jump component and they are not Lipschitz continuous. 
Moreover, these generators are not defined on the whole space $\R\times\R\times\R$.  
Using a decomposition approach based on Proposition \ref{propjeulin}, we obtain the following result whose proof is detailed in Section \ref{preuve theoreme}. 
\begin{Theorem}\label{THEUBSDEs}
The BSDEs \reff{EDSRf^1}, \reff{EDSRf^2} and \reff{EDSRf^3} admit solutions $(Y,Z,U)$, $(\Yc,\Zc,\Uc)$ and $(\Upsilon,\Xi,\Theta)$ in $\Sc^\infty_\G\times L^2_\G\times L^2(\lambda)$. Moreover $Y\in  \Sc^{\infty,+}_\G$.
\end{Theorem}

\subsection{A verification Theorem}
We now turn to the sufficient condition of optimality. 
As explained in Subsection \ref{martingalprincopt}, a candidate to be an optimal strategy is a process $\pi^*\in\Ac$ such that $J^{\pi^*}$ is a martingale, which implies that $dK^{\pi^*}=0$. This leads to 
\beqs
\pi^*_t  & = & \text{arg}\min_{\pi\in\R}  K_t(\pi) \;,
\enqs
which gives the implicit equation in $\pi^*$
\beqs
\pi_t^* & = & \big(\Yc_{t^-} - V^{x,\pi^*}_{t^-}\big)\frac{\mu_{t} Y_{t^-}+\sigma_{t} Z_{t}+\lambda^\G_{t}\beta_{t} U_{t}}{ |\sigma_{t}|^2 Y_{t^-} + \lambda^\G_{t} |\beta_{t}|^2(U_{t}+Y_{t^-})}
 + \frac{\sigma_{t}Y_{t^-}\Zc_{t} + \lambda^\G_{t}\beta_{t}\Uc_{t}(Y_{t^-}+U_{t})}{ |\sigma_{t}|^2 Y_{t^-} + \lambda^\G_{t}|\beta_{t}|^2(U_{t}+Y_{t^-})}\;.
\enqs
Integrating each side of this equality w.r.t. ${dS_t\over S_{t^-}}$ leads to the following SDE
\beq\label{EDSportfeuilleOptimal}
V^*_t & = & x + \int_0^t\big(\Yc_{r^-} - V^*_{r^-}\big)\frac{\mu_{r} Y_{r^-}+\sigma_{r} Z_{r}+\lambda^\G_{r}\beta_{r} U_{r}}{ |\sigma_{r}|^2 Y_{r^-} + \lambda^\G_{r} |\beta_{r}|^2(U_{r}+Y_{r^-})}\frac{d S_{r}}{S_{r^-}}\\
 & & +
\int_0^t
\frac{\sigma_{r}Y_{r^-}\Zc_{r} + \lambda^\G_{r}\beta_{r}\Uc_{r}(Y_{r^-}+U_{r})}{ |\sigma_{r}|^2 Y_{r^-} + \lambda^\G_{r}|\beta_{r}|^2(U_{r}+Y_{r^-})}\frac{d S_{r}}{S_{r^-}}\;, \quad t\in [0,T\wedge\tau]\;. \nonumber
\enq
We first study the existence of a solution to SDE \reff{EDSportfeuilleOptimal}. 
\begin{Proposition}\label{PropexistssoloptSDE}
 The SDE \reff{EDSportfeuilleOptimal} admits a solution $V^*$ which satisfies
\beq\label{cond int V*}
 \E\Big[\sup_{t\in [0,T\wedge \tau]} |V_t^*|^2 \Big] & < & \infty\;.
\enq
\end{Proposition}
\ni\textbf{Proof.}
To alleviate the notation we rewrite \reff{EDSportfeuilleOptimal} under the form
\begin{equation}\label{EDSV*EF}
\left\{\begin{array}{rcl}
V^*_0 & = & x  \;,\\
dV^*_t & = & (E_t V^*_{t^-} - F_t)(\mu_t dt + \sigma_t dW_t + \beta_t dM_t)\;,
\end{array}\right.
\end{equation}
where $E$ and $F$ are defined by
\beqs
E_t & := & -\frac{\mu_t Y_{t^-} + \sigma_t Z_t + \lambda^\G_t \beta_t U_t}{|\sigma_t|^2 Y_{t^-} + \lambda^\G_t |\beta_t|^2 (U_t + Y_{t^-})} \;,\\\\
F_t & := & -\frac{\lambda^\G_t \beta_t \Uc_t (Y_{t^-} + U_t) + \mu_t Y_{t^-}\Yc_{t^-} + \lambda^\G_t \beta_t U_t \Yc_{t^-} + \sigma_t Z_t \Yc_{t^-} + \sigma_t \Zc_t Y_{t^-}}{|\sigma_t|^2 Y_{t^-} + \lambda^\G_t |\beta_t|^2 (U_t + Y_{t^-})} \;,
\enqs
for all $t\in[0,T]$.
We first notice that from \textbf{(H$S$)} (ii), and since $Y\in \Sc^{\infty,+}_\G$ and $\lambda^\G (Y+U)$ is nonnegative, 
 there exists a constant $C>0$ such that 
\beqs
|\sigma_t|^2 Y_t + \lambda^\G_t |\beta_t|^2 (U_t + Y_{t^-}) & \geq & C\;,\quad \P\otimes dt-a.e.
\enqs
Therefore, using $(Y,Z,U)$, $(\Yc,\Zc,\Uc)$, $(\Upsilon, \Xi,\Theta)\in \Sc^\infty_\G\times L^2_\G\times L^2(\lambda)$, Remark \ref{remUbdd} and \textbf{(H$S$)}, we get that $E$ and $F$ are square integrable
\beqs
\E\Big[\int_0^T \Big(|E_t|^2 + |F_t|^2 \Big) dt\Big] & < & \infty\;.
\enqs
Using It\^o's formula, we obtain that the process $V^*$ defined by 
\beq\label{def-sol-V-*}
V^*_t   &   :=  &  (x + \Psi_t)\Phi_t  \;,\quad t\in[0,T\wedge \tau) \;,\\ \nonumber
  \mbox{ and }~ V^*_{T\wedge\tau} & = & \mathds{1}_{\tau\leq T}\big[(1+E_{\tau}\beta_\tau)V^*_{\tau^-}-F_{\tau}\beta_\tau\big]+ \mathds{1}_{\tau> T}(x + \Psi_T)\Phi_T  \;,
 \enq
where
\beqs
\Phi_t & := & \exp \Big( \int_0^t \Big(E_s ( \mu_s - \lambda^\G_s \beta_s)  - \frac{1}{2}  |\sigma_sE_s|^2\Big) ds + \int_0^t\sigma_s E_s dW_s \Big)\;, 
\enqs
and
\beqs
\Psi_t & := & -\int_0^t \frac{F_s}{\Phi_s} \Big[ \mu_s - \lambda^\G_s \beta_s - |E_s \sigma_s|^2 \Big]ds - \int_0^t \frac{F_s}{\Phi_s}\sigma_s dW_s\;, 
\enqs
for all $t\in[0,T]$, is solution to \reff{EDSportfeuilleOptimal}. 

\ni We now prove that $V^*$ defined by \reff{def-sol-V-*} satisfies \reff{cond int V*}. We proceed in two steps.

 \vspace{2mm}

\ni\textbf{Step 1:} We prove that 
\beq\label{condtermL2}
\E\Big[|V^*_{ T \wedge \tau}|^2\Big] & < & \infty\;.
\enq
Since $V^*$ satisfies  \reff{EDSV*EF}, we have $V^*= V^{x,\pi^*}$ where $\pi^*$ is given by
\beqs
\pi^*_t & = & E_tV^*_{t^-} - F_t\;,\quad t\in [0,T]\;.
\enqs
 We therefore have $Y|V^*_{. \wedge \tau}-\Yc|^2=J^{\pi^*}-\Upsilon$ and from \reff{decompJpi} and the dynamics of $\Upsilon$ given by \reff{EDSRf^3}, we have
 \beqs
 d \big(Y_t|V^*_{t\wedge \tau}-\Yc_t|^2\big)
 & = & dM^{*}_t+dK^{\pi^*}_t-\mathfrak{h}(t)dt
 \enqs
where $M^*$ is a locally square integrable martingale with $M^*_0=0$. From the definition of $K^{\pi^*}$  and using the fact that 
\beqs
\pi^*_t & = & \frac{X_{t^-}^{\pi^*}(\mu_tY_{t^-}+\sigma_tZ_t+\lambda^\G_t\beta_tU_t) + \sigma_tY_{t^-}\Zc_t + \lambda^\G_t\beta_t\Uc_t(Y_{t^-}+U_t)}{ |\sigma_t|^2 Y_{t^-} + \lambda^\G_t |\beta_t|^2(U_t+Y_{t^-})}\;,\quad
\enqs 
 we get $K_t^{\pi^*}=0$ for all $t\in [0,T\wedge\tau]$. Therefore, from the definition of $\mathfrak{h}$, we get
\beqs
Y_{T \wedge \tau}|V^{*}_{ T \wedge \tau} - \Yc_{T \wedge \tau}|^2 & = & Y_0 |x - \Yc_0|^2 + {M}^*_{ T \wedge \tau}   + \int_0^{T \wedge \tau} \Big[ | \Zc_t |^2 Y_t + \lambda^\G_t (U_t + Y_t) |\Uc_t|^2\\
 & & - \; \frac{ | \sigma_t Y_t \Zc_t + \lambda^\G_t \beta_t \Uc_t (U_t + Y_t)|^2}{|\sigma_t|^2 Y_t + \lambda^\G_t |\beta_t|^2(U_t +Y_t)} \Big] dt\;.
\enqs 

Since $M^*$ is a local martingale, there exists an increasing sequence of $\G$-stopping times $(\nu_i)_{i \in \N}$ such that $\nu_i \rightarrow +\infty$ as $i\rightarrow \infty$ and
\beq \label{egalite auxiliaire}
\nonumber \E \big[Y_{{T\wedge \tau} \wedge \nu_i} |V^{*}_{{T \wedge \tau} \wedge \nu_i} - \Yc_{{T\wedge \tau}\wedge \nu_i}|^2 \big] & = & Y_0 |x - \Yc_0|^2 + \E   \int_0^{{T\wedge \tau} \wedge \nu_i} \Big[ | \Zc_t |^2 Y_t + \lambda^\G_t (U_t + Y_t) |\Uc_t|^2 
\\
&& - \frac{ | \sigma_t Y_t \Zc_t + \lambda^\G_t \beta_t \Uc_t (U_t + Y_t)|^2}{|\sigma_t|^2 Y_t + \lambda^\G_t |\beta_t|^2(U_t +Y_t)} \Big] dt \;.
\enq
Since $Y\in  \Sc^{\infty,+}_\G$, there exists a positive constant $C$ such that 
\beqs
 \E \big[|V^{*}_{{T \wedge \tau} \wedge \nu_i} - \Yc_{{T \wedge \tau} \wedge \nu_i}|^2 \big] & \leq & C \E \big[Y_{{T\wedge \tau} \wedge \nu_i}|V^{*}_{{T \wedge \tau} \wedge \nu_i} - \Yc_{{T \wedge \tau} \wedge \nu_i}|^2 \big] \;.
 \enqs
Therefore, using \reff{egalite auxiliaire}, we get that 
\beqs
 \E \big[|V^{*}_{{T \wedge \tau} \wedge \nu_i} - \Yc_{{T \wedge \tau} \wedge \nu_i}|^2 \big] & \leq & C\Big(Y_0 |x - \Yc_0|^2 + \E   \int_0^T \Big[ | \Zc_t |^2 Y_t + \lambda^\G_t (U_t + Y_t) |\Uc_t|^2 \Big] dt \Big) \;.
\enqs
Since $Y$, $U$ and $\Uc$ are uniformly bounded and $\Zc \in L^2_\G$, there exists a constant $C$ such that
\beqs
\E \big[\big|V^{*}_{{T\wedge \tau} \wedge \nu_i} - \Yc_{{T\wedge \tau} \wedge \nu_i}\big|^2 \big] & \leq & C\;.
\enqs
From Fatou's lemma, we get that 
\beqs
\E \big[ \big|V^{*}_{T\wedge \tau} - \Yc_{T\wedge \tau}\big|^2 \big] & \leq & \lim_{i \;\rightarrow} \inf_{\hspace{-2mm}\infty}\E \big[\big|V^{*}_{{T\wedge \tau} \wedge \nu_i} - \Yc_{{T\wedge \tau} \wedge \nu_i}\big|^2 \big] ~\leq ~ C \;.
\enqs
Which implies that
\beqs
\E\big[\big|V^{*}_{T\wedge \tau}\big|^2\big] & \leq & C + 2 \E \big[V^{*}_{T\wedge \tau}\Yc_{T\wedge \tau}\big]\;.
\enqs
Finally, using the Young inequality and noting that $\Yc$ is uniformly bounded, it follows that there exists a constant $C$ such that
\beqs
\E\big[\big|V^{*}_{T\wedge \tau}\big|^2\big] & \leq & C \;.
\enqs
\textbf{Step 2:} We prove that 
\beqs
\E\Big[ \sup_{t \in [0, T\wedge \tau]} |V^{*}_{t}|^2 \Big]& < & \infty \;.
\enqs
For that we remark that $V^*_{.\wedge\tau}$ is solution to the following linear BSDE
\beq\label{BSDE V*}
V^*_{t\wedge\tau} & = & V^*_{T\wedge\tau} -\int_{t\wedge\tau}^{T\wedge\tau}{\mu_s\over \sigma_s} z_s ds -\int_{t\wedge\tau}^{T\wedge\tau} z_s  dW_s - \int_{t\wedge\tau}^{T\wedge\tau}u_s dM_s\;,\quad t\in[0,T]\;,
\enq
with 
\beqs
z_t & := & \sigma_t \frac{(\Yc_{t^-} - V^*_{t^-})(\mu_tY_{t^-}+\sigma_tZ_t+\lambda^\G_t\beta_tU_t) + \sigma_tY_{t^-}\Zc_t + \lambda^\G_t\beta_t\Uc_t(Y_{t^-}+U_t)}{ |\sigma_t|^2 Y_{t^-} + \lambda^\G_t |\beta_t|^2(U_t+Y_{t^-})}\;,\\\\
u_t & := & \beta_t \frac{(\Yc_{t^-} - V^*_{t^-})(\mu_t Y_{t^-}+\sigma_tZ_t+\lambda^\G_t\beta_tU_t) + \sigma_tY_{t^-}\Zc_t + \lambda^\G_t\beta_t\Uc_t(Y_{t^-}+U_t)}{ |\sigma_t|^2 Y_{t^-} + \lambda^\G_t |\beta_t|^2(U_t+Y_{t^-})}\;,
\enqs
for all $t\in[0,T]$.
Applying It\^o's formula to $|V^*|^2$ we have
\beqs
\E|V^*_{t\wedge \tau}|^2 & = & \E|V^*_{T\wedge \tau}|^2-2\E\int_{t\wedge\tau}^{T\wedge\tau}V^*_{s\wedge \tau}{\mu_s\over \sigma_s} z_s  ds 
  -\E\int_{t\wedge\tau}^{T\wedge\tau}|z_s|^2ds-\E\int_{t\wedge\tau}^{T\wedge\tau}|u_s|^2\lambda_sds\;,
\enqs
for all $t\in [0,T]$. Using \reff{condtermL2}, \textbf{(H$S$)} and the  Young inequality we obtain the existence of a constant $C$ such that
\beqs
\E|V^*_{t\wedge \tau}|^2  +\E\int_{t\wedge\tau}^{T\wedge\tau}|z_s|^2ds+\E\int_{t\wedge\tau}^{T\wedge\tau}|u_s|^2\lambda_sds& \leq & C\Big(1+ \E\int_{t}^{T}|V^*_{s\wedge \tau}|^2\Big)\;.
\enqs
We then deduce from the Gronwall inequality that
\beq\label{ineg inter V*S2}
\sup_{t\in[0,T]}\E|V^*_{t\wedge \tau}|^2  +\E\int_{0}^{T\wedge\tau}|z_s|^2ds+\E\int_{0}^{T\wedge\tau}|u_s|^2\lambda_sds & < & +\infty\;.
\enq
Now from \reff{BSDE V*}, we have
\beqs
\E\Big[\sup_{t\in[0,T]}|V^*_{t\wedge \tau}|^2\Big] & \leq & 3\Big(|V^*_{0}|^2+\E\Big[\sup_{t\in[0,T]}\Big|\int_{0}^{t\wedge\tau}{\mu_s\over \sigma_s}z_sds\Big|^2\Big]\\
 & & +\E\Big[\sup_{t\in[0,T]}\Big|\int_0^{t\wedge \tau}z_sdW_s+\int_0^{t\wedge \tau}u_sdM_s\Big|^2\Big]\Big)\;.
\enqs
From \textbf{(H$S$)} and the  BDG inequality, there exists a constant $C$ such that
\beqs
\E\Big[\sup_{t\in[0,T]}|V^*_{t\wedge \tau}|^2\Big] & \leq & C\Big(1 +\E\int_{0}^{T\wedge\tau}|z_s|^2ds+\E\int_{0}^{T\wedge\tau}|u_s|^2\lambda_sds\Big)\;.
\enqs
This last inequality with \reff{ineg inter V*S2} gives \reff{cond int V*}.
\ep

\vspace{2mm}

\ni As explained previously, we now consider the strategy $\pi^*$ defined by  
\beq\label{strategie optimale}
\pi^*_t & = & \frac{(\Yc_{t^-} - V^*_{t^-})(\mu_tY_{t^-}+\sigma_tZ_t+\lambda^\G_t\beta_tU_t) + \sigma_tY_{t^-}\Zc_t + \lambda^\G_t\beta_t\Uc_t(Y_{t^-}+U_t)}{ |\sigma_t|^2 Y_{t^-} + \lambda^\G_t |\beta_t|^2(U_t+Y_{t^-})}\;,\quad
\enq
for all $t\in[0,T]$. We first notice from the expressions of $\pi^*$ and $V^*$ that 
\beq\label{indentificationV*etVpi*}
V^{x,\pi^*}_t & = & V^*_t \;,
\enq
for all $t\in[0,T]$. Using \reff{cond int V*} and \reff{indentificationV*etVpi*}, we have 
\beq\label{cond int Vpi*}
 \E\Big[\sup_{t\in [0,T\wedge \tau]} |V_t^{x,\pi^*}|^2 \Big] & < & \infty\;.
\enq
We can now state our verification theorem which is the main result of this section.  
\begin{Theorem}\label{Th Verif} 
 The strategy $\pi^*$ given by \reff{strategie optimale} belongs to the set $\Ac$ and is optimal for the mean-variance problem \reff{MVpb}. Thus we have
\beqs
\E \Big[\big|V_{T\wedge\tau}^{x,\pi^*}-H\big|^2 \Big] & = & \min_{\pi\in\Ac}\E \Big[\big|V_{T\wedge\tau}^{x,\pi}-H\big|^2 \Big]~=~ Y_0|x-\Yc_0|^2+\Upsilon_0\;,
\enqs
where $Y,\Yc$ and $\Upsilon$ are solutions to \reff{EDSRf^1}-\reff{EDSRf^2}-\reff{EDSRf^3}.
\end{Theorem}
\ni To prove this verification theorem, we first need of the following lemma.
\begin{Lemma}\label{lem M pi ml}
For any $\pi\in\Ac$, the process $M^\pi_{.\wedge\tau}$ defined by \reff{decompJpi} is a $\G$-local martingale.
\end{Lemma}
\ni\textbf{Proof.} Fix $\pi\in \Ac$. Then from the definition of $V^{x,\pi}$, \textbf{(H$S$)} and the BDG inequality, we have 
\beq\label{estim V S2}
\E\Big[\sup_{t\in[0,T]}\big|V^{x,\pi}_{t\wedge\tau}\big|^2\Big] & < & \infty\;.
\enq
Define the sequence of $\G$-stopping times $(\nu_n)_{n \geq 1}$ by 
\beqs
\nu_n & := & \inf\Big\{ s\geq 0~:~ \big|V^{x,\pi}_{s\wedge\tau}\big|\geq n \Big\}\;,
\enqs
for all $n\geq 1$. First, notice that $(\nu_n)_{n \geq 1}$ is nondecreasing and goes to infinity as $n$ goes to infinity from \reff{estim V S2}. Moreover, 
from the definition of $\nu_n$, 
we have 
\beqs
|V_s^{x,\pi}\mathds{1}_{s\in [0,\nu_n\wedge \tau)}| & \leq & n 
\enqs
for all $s\in[0,T]$. 
 Then, since $\pi\in\Ac$, $Y,\Yc\in \Sc^\infty_\G$ and $Z, \Zc, \Xi\in L^2_\G$, we get 
\beqs
\E\Big[\int_0^{\tau\wedge\nu_n\wedge T}\Big|2Y_t X^\pi_t ( \pi_t \sigma_t - \Zc_t ) + Z_t |X^\pi_t|^2 + \Xi_t\Big|^2dt \Big] & < & \infty\;,
\enqs
for all $n\geq 1$. Moreover, since $U, \Uc, \Theta\in L^2(\lambda)$, we get from Remark \ref{remUbdd}  
\beqs
\E\Big[\int_0^{\tau\wedge\nu_n\wedge T}\Big|(2 X^\pi_{t^-}+ \pi_t \beta_t - \Uc_t ) ( \pi_t \beta_t - \Uc_t ) ( Y_{t^-} + U_t)  
+{|X^\pi_{t^-}|^2U_t}+ \Theta_t\Big|\lambda^\G_tdt \Big] & < & \infty\;,
\enqs
for all $n\geq 1$. Therefore, we get that the stopped process $M^\pi_{.\wedge\tau\wedge\nu_n}$ is a $\G$-martingale. 
\ep

\vspace{2mm}

\ni\textbf{Proof of Theorem \ref{Th Verif}.} As explained in Subsection \ref{subsectionduale},  we check each of the points (i), (ii), (iii) and (iv).

\ni (i)  From the definition of $Y$, $\Yc$ and $\Upsilon$, we have 
\beqs
J^{\pi}_{T\wedge\tau} & = & Y_{T\wedge\tau} \big|V_{T\wedge\tau}^{x,\pi}-H \big|^2+\Upsilon_{T\wedge\tau}  ~ = ~ {\big| V_{T\wedge\tau}^{x,\pi} -H\big|}^2\;,
\enqs
for all $\pi\in \Ac$.

\vspace{2mm}

\ni (ii)  From the definition of the family $(J^\pi)_{\pi\in\Ac}$, we have
\beqs
J^\pi_0 & = & Y_0|V_0^{x,\pi}-\Yc_0|^2+\Upsilon_0~ = ~ Y_0|x-\Yc_0|^2+\Upsilon_0\;,
\enqs
for all $\pi\in \Ac$.

\vspace{2mm}

\ni (iii) Fix $\pi\in \Ac$. 
Since $Y,\Yc,\Upsilon\in  \Sc^\infty_\G$, we have from the definition of $J^\pi$ and the BDG inequality 
\beq\label{Jpi de classe D}
\E\Big[\sup_{t\in[0,T]}|J^\pi_t|\Big] & < & +\infty\;.
\enq
Now, fix $s,t\in [0,T]$ such that $s\leq t$. Using the decomposition \reff{decompJpi} and Lemma \ref{lem M pi ml}, there exists an increasing sequence of $\G$-stopping times $(\nu_i)_{i\geq 1}$ such that $\nu_i\rightarrow+\infty$ as $i\rightarrow+\infty$ and 
\beq\label{ppte smgjpi}
\E\Big[J^\pi_{t\wedge\nu_i}\big|\Gc_s\Big] & \geq & J^\pi_{s\wedge\nu_i}\;,
\enq
for all $i\geq 1$. Then, from \reff{Jpi de classe D}, we can apply the conditional dominated convergence theorem and we get by sending $i$ to $\infty$ in \reff{ppte smgjpi}
\beqs
\E\Big[J^\pi_{t}\big|\Gc_s\Big] & \geq & J^\pi_{s}\;,
\enqs
for all $s,t\in[0,T]$ with $s\leq t$. 

\vspace{2mm}

\ni (iv) We now check that $\pi^*\in\Ac$  
i.e. $\E\int_0^{T \wedge \tau}|\pi^*_s|^2ds<\infty$. Using the definition of $\pi^*$ and \reff{indentificationV*etVpi*} we have  that 
$V^{x,\pi^*}$ is solution to the linear BSDE
\beqs
V^{x,\pi*}_t & = & V^{x,\pi^*}_{T\wedge\tau} -\int_{t\wedge\tau}^{T\wedge\tau}{\mu_s\over \sigma_s} z_s ds -\int_{t\wedge\tau}^{T\wedge\tau} z_s  dW_s - \int_{t\wedge\tau}^{T\wedge\tau}u_s dM_s\;,\quad t\in[0,T]\;,
\enqs
with 
\beqs
z_t ~ = ~ \sigma_t \pi^*_t
 & \mbox{ and } & u_t ~ = ~ \beta_t  \pi^*_t\;,
 \enqs
for all $t\in[0,T]$.
Therefore, using \reff{cond int Vpi*}, \textbf{(H$S$)}, applying It\^o's formula to $|V^{x,\pi*}|^2$, using the Young inequality, the BDG inequality and the Gronwall inequality (see e.g. the proof of Proposition 2.2 in \cite{BBP97}), we get 
\beqs
\E\Big[\int_0^{T\wedge\tau}|\pi^*_s|^2ds\Big] & < & \infty\;.
\enqs

We now check that $J^{\pi^*}$ is a $\G$-martingale.  
Since $K^{\pi^*}$ is constant, we obtain from Lemma \ref{lem M pi ml} that $J^{\pi^*}$ is a  $\G$-local martingale. 
Then, from the expression of $J^{\pi^*}$ and since $Y,\Yc,\Upsilon\in  \Sc^\infty_\G$, there exists a constant $C$ such that
\beqs
\E\Big[\sup_{t\in[0,T]}|J_t^{\pi^*}|\Big] & \leq & C\Big(1+\E\Big[\sup_{t\in[0,T\wedge\tau]}|V_t^{x,\pi^*}|^2\Big]\Big)\;.
\enqs
Using \reff{cond int Vpi*}, 
we get that 
\beqs
\E\Big[\sup_{t\in[0,T]}|J_t^{\pi^*}|\Big] & < & +\infty\;.
\enqs
Therefore, $J^{\pi^*}$ is a true $\G$-martingale and $\pi^*$ is optimal. \ep

\section{A decomposition approach for solving BSDEs in the filtration $\G$}
\label{preuve theoreme}
We now prove Theorem \ref{THEUBSDEs} via a decomposition procedure. We first provide a general result which gives existence of a solution to a BSDE in the enlarged filtration $\G$ as soon as an associated BSDE in the filtration $\F$ admits a solution. Actually the associated BSDE is defined by the terms appearing in the decomposition of the coefficients of the BSDE in $\G$ given by Proposition \ref{propjeulin}.
 We therefore introduce the spaces of processes where solutions in $\F$ classically lie.
\begin{itemize}
\item $\Sc_{\F}^\infty$ is the subset of $\R$-valued  continuous $\F$-adapted processes $(Y_t)_{t\in[0,T]}$ essentially bounded
\beqs
{\| Y \|}_{\Sc^\infty} & := & \Big\|\sup_{t\in[0,T]}|Y_{t}|\Big\|_\infty~<~\infty \;.
\enqs
\item  $\Sc_{\F}^{\infty,+}$ is the subset of  $\Sc_{\F}^\infty$ of processes $(Y_t)_{t\in[0,T]}$ valued in $(0,\infty)$, such that
\beqs
{\Big\| {1\over Y} \Big\|}_{\Sc^\infty} & < & \infty \;.
\enqs
\item $L^2_{\F}$ is the subset of $\R$-valued  $\Pc(\F)$-measurable processes $(Z_t)_{t\in[0,T]}$ such that
\beqs
\|Z\|_{L^2}  & := & \Big(\E\Big[ \int_0^T |Z_t|^2 dt \Big]\Big)^{1\over2} ~< ~\infty \;.
\enqs
\end{itemize}
Finally since the BSDEs associated to our mean-variance problem have generators with superlinear growth, we consider the additional space of BMO-martingales: $\mathrm{BMO}(\P)$ is the subset of $(\P,\F)$-martingales $m$ such that
\beqs
 \|m\|_{\mathrm{BMO}(\P)}  & := & \sup_{\nu\in \Tc_\F[0,T]} \Big\|\E \big[\langle m \rangle_T - \langle m \rangle_\nu | \Fc_\nu\big]^{1\over2} \Big\|_\infty~<~ \infty \;,
\enqs
where $\Tc_\F[0,T]$ is the set of $\F$-stopping times on $[0,T]$. 
 This means local martingales of the form $m_t = \int_0^t Z_s dW_s$, $t\in[0,T]$, are $\rm{BMO}(\P)$-martingale if and only if
\beqs
 \Big\|\int_0^.Z_sdW_s\Big\|_{\mathrm{BMO}(\P)}  & := & \sup_{\nu\in \Tc_\F[0,T]}{\Big\|\Big(\E\Big[ \int_\nu^T |Z_t|^2 dt\Big|\Fc_\nu \Big]\Big)^{1\over2}\Big\|}_\infty ~< ~\infty \;.
\enqs
\subsection{A general existence theorem for BSDEs with random horizon}
We provide here a general result on existence of a solution to a BSDE driven by $W$ and $N$ with horizon $T\wedge \tau$. 
We consider a generator function $F:\Omega\times [0,T]\times \R\times\R \times \R\rightarrow \R$, which is $\Pc(\G)\otimes \Bc(\R)\otimes\Bc(\R)\otimes \Bc(\R)$-measurable, and a terminal condition $\xi$ which is a $\Gc_{T\wedge \tau}$-measurable random variable of the form 
\beq\label{decomxi}
\xi & = & \xi^b\mathds{1}_{T< \tau}+ \xi^a_\tau\mathds{1}_{T\geq \tau}\;,
\enq
where $\xi^b$ is an $\Fc_T$-measurable bounded  random variable and $\xi^a\in \Sc^\infty_\F$.   From Proposition \ref{propjeulin} and Remark \ref{remjeulin}, we can write
\beq\label{decomposition f}
F(t,.)\mathds{1}_{t\leq \tau} & = & F^b(t,.)\mathds{1}_{t\leq \tau}\;,\quad t\geq 0\;, 
\enq
where $F^b$ is a $\Pc(\F)\otimes\Bc(\R)\otimes\Bc(\R)\otimes\Bc(\R)$-measurable map. We then introduce the following BSDE
\beq\label{EDSR-avant-saut}
Y^b_t & = & \xi^b+\int_t^TF^b(s,Y^b_s,Z^b_s,\xi^a_s-Y^b_s)ds-\int_t^TZ_s^bdW_s\;,\quad t\in[0,T]\;.\enq
\begin{Theorem}\label{Thm-exist-gen}
Assume that the BSDE \reff{EDSR-avant-saut} admits a solution $(Y^b,Z^b)\in \Sc^\infty_\F\times L^2_\F$. Then BSDE
\beq\label{EDSRG}
\hspace{-9mm} Y_t & = & \xi + \int_{t\wedge \tau}^{T\wedge \tau}F(s,Y_s,Z_s,U_s)ds  -\int_{t\wedge \tau}^{T\wedge \tau}Z_sdW_s-\int_{t\wedge \tau}^{T\wedge \tau}U_s dN_s\;, \quad t\in [0,T]\;,
\enq
admits a solution $(Y,Z,U)\in  \Sc^\infty_\G\times L^2_\G\times L^2(\lambda) $ given by 
\beq \nonumber
Y_t & = & Y_t^b\mathds{1}_{t<\tau}+ \xi^a_\tau\mathds{1}_{t\geq\tau}\;,\\\label{defsolEDSRG}
Z_t & = & Z^b_t\mathds{1}_{t \leq \tau}\;,\\
U_t & = &\big( \xi^a_t-Y^b_t\big)\mathds{1}_{t \leq \tau}\;,\nonumber
\enq
for all $t\in [0,T]$. 
\end{Theorem}
\ni \textbf{Proof.}
We proceed in three steps. 

\noindent\textbf{Step 1:} We prove that for $t$ $\in$ $[0,T]$, $(Y, Z, U)$ defined by \reff{defsolEDSRG} satisfies the equation \reff{EDSRG}.   
We distinguish  three cases. \\
\textbf{Case 1:} $\tau>T$.

\ni From \reff{defsolEDSRG}, we get $Y_t=Y^b_t$, $Z_t=Z^b_t$ and $U_t=\xi^a_t-Y^b_t$ for all $t\in [0,T]$.  Then, using that $(Y^b,Z^b)$ is a solution to \reff{EDSR-avant-saut}, we have 
\beqs
Y_t & = & \xi^b + \int_t^T F^b(s,Y_s, Z_s, U_s)ds-\int_t^T Z^b_sdW_s \;.
\enqs
Since the predictable processes $Z$ and $Z^b$ are indistinguishable on $\{\tau> T\}$, we have from Theorem 12.23 of \cite{HeW},  $\int_t^TZ_sdW_s = \int_t^TZ_s^bdW_s $ on $\{\tau > T\}$.
Moreover since $\xi = \xi^b$ and $\int_{t\wedge \tau}^{T\wedge\tau} U_sdN_s=0$ on $\{\tau>T\}$ we get by using \reff{decomposition f}
\beqs
Y_t& = &\xi+\int_{t\wedge \tau}^{T\wedge \tau}F(s,Y_s,Z_s,U_s)ds-\int_{t\wedge \tau}^{T\wedge \tau}Z_sdW_s-\int_{t\wedge \tau}^{T\wedge \tau}U_sdN_s\;.
\enqs

\vspace{2mm}

\ni\textbf{Case 2:} $\tau\in (t,T]$. 

\ni From \reff{defsolEDSRG}, we have $Y_t =Y_t^b$. Since $(Y^b,Z^b)$ is solution to \reff{EDSR-avant-saut}, we have 
\beqs
Y_t & = & {Y}^b_\tau+\int_t^\tau F^b(s,Y_s^b, Z_s^b, \xi^a_s-Y_s^b)ds-\int_t^\tau Z^b_sdW_s \;.
\enqs
 Still using \reff{decomposition f} and \reff{defsolEDSRG}, we get
\beqs
Y_t & = & \xi^a_\tau +\int_t^\tau F(s,Y_s, Z_s, U_s)ds-\int_t^\tau Z^b_sdW_s -(\xi^a_\tau-Y^b_\tau)\;.
\enqs
Since the predictable processes $Z\mathds{1}_{. <\tau}$ and $Z^b\mathds{1}_{. <\tau}$ are indistinguishable on $\{\tau> t\}\cap\{\tau\leq T\}$, we have from Theorem 12.23 of \cite{HeW},  $\int_t^{T\wedge\tau}Z_sdW_s = \int_t^{T\wedge\tau}Z_s^bdW_s $ on $\{\tau> t\}\cap\{\tau\leq T\}$. Therefore, we get
\beqs
Y_t & = & \xi^a_\tau+\int_t^\tau F(s,Y_s, Z_s, U_s)ds-\int_t^\tau Z_sdW_s -(\xi^a_\tau-Y^b_\tau)\;.
\enqs
 Finally, we easily check from the definition of $U$ that $\int_t^{T \wedge \tau}U_sdN_s=\xi^a_\tau-Y^b_\tau$. Therefore, we get using   \reff{decomxi}
\beqs
 Y_t & = & \xi +\int_{t\wedge \tau}^{T\wedge \tau } F(s,Y_s, Z_s, U_s)ds-\int_{t\wedge \tau}^{T\wedge \tau } Z_sdW_s - \int_{t\wedge \tau}^{T\wedge \tau}U_sdN_s\;.
\enqs
 
\vspace{2mm}

\ni\textbf{Case 3:} $\tau\leq t$. 

\ni Then, from \reff{defsolEDSRG}, we have $Y_t=\xi^a_\tau$. We therefore get on $\{ \tau\leq t\}$ by using \reff{decomxi}
\beqs
Y_t& = &\xi+\int_{t\wedge \tau}^{T\wedge \tau}F(s,Y_s,Z_s,U_s)ds-\int_{t\wedge \tau}^{T\wedge \tau}Z_sdW_s-\int_{t\wedge \tau}^{T\wedge \tau}U_sdN_s\;.
\enqs

\vspace{2mm}

\noindent\textbf{Step 2:} We notice that $Y$ is a c\`ad-l\`ag $\G$-adapted process and $U$ is $\Pc(\G)$-measurable since $Y^b$ and $\xi^a$ are continuous and $\F$-adapted. 
We also notice  from its definition that the process  $Z$ is $\Pc(\G)$-measurable,  since 
 $Z^b$ is $\Pc({\F})$-measurable.

\vspace{2mm}

\noindent\textbf{Step 3:} We now prove that the solution satisfies the integrability conditions. 
\begin{itemize}
\item From the definition of $Y$, we have 
\beq\label{maj Y}
|Y_t|  & \leq & |Y_t^b|+ |\xi^a_t| \;,\quad t\in [0,T]\;.
\enq
Since $Y^b\in \Sc^\infty_\F$ and $\xi^a\in \Sc^\infty_\F$, 
we get that $\|Y\|_{\Sc^\infty}<+\infty$.
\item From the definition of the process $Z$, we have $Z \in L^2_\G$.

\item From the definition of $U$, we have 
\beqs
|U_t|  & \leq & |Y_t^b|+ |\xi^a_t|\;,\quad t\in [0,T]\;.
\enqs
Since $Y^b\in \Sc^\infty_\F$, $\xi^a\in \Sc^\infty_\F$ and $\lambda$ is bounded, we get $U \in L^2(\lambda)$.
\end{itemize}

\ep

\vspace{2mm}

\ni Using this abstract result we prove the existence of solutions to each of the BSDEs  \reff{EDSRf^1}, \reff{EDSRf^2} and \reff{EDSRf^3}
in the following subsections. 

\subsection{Solution to the BSDE $(\mathfrak{f}, 1)$}
 Following Theorem \ref{Thm-exist-gen}, we consider for coefficients $(\mathfrak{f}, 1)$ the BSDE in $\F$: find $(Y^b,Z^b)\in \Sc^\infty_\F\times L^2_\F$ such that
\begin{equation}\label{EDSR0}
\left\{\begin{array}{rcl}
d Y^b_t & = & \Big\{ \displaystyle{\frac{|(\mu_t - \lambda_t \beta_t) Y^b_t + \sigma_t Z^b_t + \lambda_t \beta_t |^2}{|\sigma_t|^2 Y^b_t + \lambda_t |\beta_t|^2}} - \lambda_t + \lambda_t Y^b_t \Big\}dt + Z^b_t d W_t\;,\quad t\in[0,T]\;,\\
Y^b_{T} & = & 1 \;.
\end{array}\right.
\end{equation}
To solve this BSDE, we have to deal with two main issues. The first is that the generator $\mathfrak{f}$ has a superlinear growth. The second difficulty is that the generator value is not defined for all the values that the process $Y$ can take. In particular the generator may explode if the process $Y$ goes to zero.  Taking in consideration these issues we get the following result. 
\begin{Proposition}\label{sol f1avsaut}
The BSDE \reff{EDSR0} has a solution $(Y^b, Z^b)$ in $\Sc^{\infty,+}_\F \times L^2_\F $ with $\int_0^.Z^bdW\in \rm{BMO}(\P)$.
\end{Proposition}

\ni\textbf{Proof.}
We first notice that the BSDE \reff{EDSR0} can be written under the form
 \begin{equation*}
\left\{\begin{array}{rcl}
d Y^b_t & = & \displaystyle{\Big\{ \frac{|\mu_t - \lambda_t \beta_t|^2}{|\sigma_t|^2} Y^b_t -\frac{\lambda_t| \beta_t|^2}{|\sigma_t|^4} |\mu_t - \lambda_t \beta_t|^2  - \lambda_t + \lambda_t Y^b_t + \frac{2(\mu_t - \lambda_t \beta_t)}{|\sigma_t|^2}(\sigma_t Z^b_t + \lambda_t \beta_t)}  \\&&\displaystyle{+\; \frac{\big|\sigma_t Z^b_t + \lambda_t \beta_t + (\lambda_t \beta_t - \mu_t) \frac{\lambda_t |\beta_t|^2}{|\sigma_t|^2}\big|^2}{|\sigma_t|^2 Y^b_t + \lambda_t |\beta_t|^2} \Big\}dt + Z^b_t d W_t}\;,\quad t\in[0,T]\;,\\
Y^b_{T} & = & 1\;.
\end{array}\right.
\end{equation*}
Since the variable $Y^b$ appears in the denominator we can not directly solve this BSDE. We then proceed in four steps. We first introduce a modified BSDE with a lower bounded denominator to ensure that the generator is well defined. We then prove via a change of probability and a comparison theorem that the solution of the modified BSDE satisfies the initial BSDE. 

\vspace{2mm} 

\ni \textbf{Step 1}: \textit{Introduction of the modified BSDE}.

Let $(Y^\eps, Z^\eps)$ be the solution in $\Sc^{ \infty}_\F \times L^2_\F $  to the BSDE
 \begin{equation}\label{EDSR modif}
\left\{\begin{array}{rcl}
d Y^\eps_t & = & \displaystyle{\Big\{ \frac{|\mu_t - \lambda_t \beta_t|^2}{|\sigma_t|^2} Y^\eps_t -\frac{\lambda_t |\beta_t|^2}{|\sigma_t|^4} |\mu_t - \lambda_t \beta_t|^2  - \lambda_t + \lambda_t Y^\eps_t + \frac{2(\mu_t - \lambda_t \beta_t)}{|\sigma_t|^2}(\sigma_t Z^\eps_t + \lambda_t \beta_t)}\\&&  \displaystyle{+ \; \frac{\big|\sigma_t Z^\eps_t + \lambda_t \beta_t + (\lambda_t \beta_t - \mu_t) \frac{\lambda_t |\beta_t|^2}{|\sigma_t|^2}\big|^2}{|\sigma_t|^2 (Y^\eps_t \vee \eps) + \lambda_t |\beta_t|^2} \Big\}dt + Z^\eps_t d W_t}\;,\quad t\in[0,T]\;,\\
Y^\eps_{T} & = & 1\;,
\end{array}\right.
\end{equation}
where $\eps$ is a positive constant such that 
\beq\label{condeps}
\exp\Big(-\int_0^T \Big( \lambda_t + \frac{|\mu_t - \lambda_t \beta_t|^2}{|\sigma_t|^2} \Big)dt\Big) & \geq & \eps\;, \quad \P-a.s.
\enq
Such a constant exists from \textbf{(H$S$)}.
Since the BSDE \reff{EDSR modif} is a quadratic BSDE, there exists a solution $(Y^\eps, Z^\eps)$ in $\Sc^{ \infty}_\F \times L^2_\F $ from \cite{kob00}. 

\vspace{2mm}

\ni \textbf{Step 2}: BMO \textit{property of the solution.}

In this part we prove that $\int_0^.Z^\eps dW \in \rm{BMO}(\P)$. Let $k$ denote the lower bound of the uniformly bounded process $Y^\eps$. Applying It\^{o}'s formula to $|Y^\eps - k|^2$, we obtain
\beq\label{estim Y eps - k}
\E \Big[ \int_\nu^T | Z^\eps_s|^2 ds \Big| \Fc_\nu \Big] & = & |1-k|^2 - | Y^\eps_\nu -k|^2 - 2 \E \Big[ \int_\nu^T (Y^\eps_s - k) f^\epsilon(s, Y^\eps_s, Z^\eps_s) ds \Big| \Fc_\nu \Big]\;,\qquad \qquad 
\enq
 for any stopping times $\nu \in \Tc_\F[0,T]$, with
\beqs
f^\epsilon(t,y,z) & = &  \frac{|\mu_t - \lambda_t \beta_t|^2}{|\sigma_t|^2} y -\frac{\lambda_t |\beta_t|^2}{|\sigma_t|^4} |\mu_t - \lambda_t \beta_t|^2  - \lambda_t + \lambda_t y + \frac{2(\mu_t - \lambda_t \beta_t)}{|\sigma_t|^2}(\sigma_t z + \lambda_t \beta_t)\\&& + \frac{|\sigma_t z + \lambda_t \beta_t + (\lambda_t \beta_t - \mu_t) \frac{\lambda_t |\beta_t|^2}{|\sigma_t|^2}|^2}{|\sigma_t|^2 (y \vee \eps) + \lambda_t |\beta_t|^2} \;,
\enqs
for all $(t,y,z)\in[0,T]\times\R\times\R$.
We can see that 
\beq\label{minoration driver eps}
f^\epsilon(t,y,z) & \geq &  I_t + G_t y + H_t z\;,
\enq
for all $(t,y,z)\in[0,T]\times\R\times\R$ where the processes $I$, $G$ and $H$ are given by
 \begin{equation*}
\left\{\begin{array}{rcl}
I_t & := & - \displaystyle{\frac{\lambda_t |\beta_t|^2}{|\sigma_t|^4} |\mu_t - \lambda_t \beta_t|^2  - \lambda_t   + 2\lambda_t \beta_t \frac{(\mu_t - \lambda_t \beta_t)}{|\sigma_t|^2}} \;, \\
G_t & := &  \displaystyle{\frac{|\mu_t - \lambda_t \beta_t|^2}{|\sigma_t|^2} + \lambda_t} \;,\\
H_t & := &   \displaystyle{2\frac{(\mu_t - \lambda_t \beta_t)}{\sigma_t}}\;,
\end{array}\right.
\end{equation*}
for all $t\in [0,T]$. We first notice that from \textbf{(H$S$)}, the processes $I$, $J$ and $K$ are bounded. 
Using  \reff{estim Y eps - k} and \reff{minoration driver eps}, we get the following inequality
\beqs
\E \Big[ \int_\nu^T | Z^\eps_s|^2 ds \Big| \Fc_\nu \Big] & \leq & |1-k|^2 - 2 \E \Big[ \int_\nu^T (Y^\eps_s - k) (I_s + G_s Y^\eps_s + H_s Z^\eps_s) ds \Big| \Fc_\nu \Big] \;.
\enqs
From the inequality $2ab\leq a^2+b^2$ for $a,b \geq0$, we get 
\beqs
\E \Big[ \int_\nu^T | Z^\eps_s|^2 ds \Big| \Fc_\nu \Big] & \leq & |1-k|^2 - 2 \E \Big[ \int_\nu^T (Y^\eps_s - k) (I_s + G_s Y^\eps_s) ds \Big| \Fc_\nu \Big] \\
& & +\; 2 \E \Big[ \int_\nu^T |H_s|^2 |Y^\eps_s - k |^2 ds \Big| \Fc_\nu \Big] + \frac{1}{2} \E \Big[ \int_\nu^T | Z^\eps_s|^2 ds \Big| \Fc_\nu \Big]\;.
\enqs
Since $I$, $G$, $H$ and $Y^\eps$ are uniformly bounded, we get
\beqs
\E \Big[ \int_\nu^T | Z^\eps_s|^2 ds \Big| \Fc_\nu \Big] & \leq & C\;,
\enqs
for some constant $C$ which does not depend on $\nu$. Therefore, $\int_0^.Z^\eps dW\in \textrm{BMO}(\P)$. 
\vspace{2mm}

\ni \textbf{Step 3}: \textit{Change of probability.}

Define the process $L^\eps$ by 
\beqs
L_t^\eps & := & 2 \frac{(\mu_t - \lambda_t \beta_t)}{\sigma_t} + 2 \frac{\sigma_t\big(\lambda_t \beta_t + \frac{\lambda_t |\beta_t|^2}{|\sigma_t|^2}(\lambda_t \beta_t - \mu_t)\big)}{|\sigma_t|^2(Y^\eps_t \vee \eps) + \lambda_t |\beta_t|^2}  + \frac{|\sigma_t|^2Z^\eps_t}{|\sigma_t|^2(Y^\eps_t \vee \eps) + \lambda_t |\beta_t|^2}  \;,
\enqs
for all $t\in[0,T]$.
Since $Y^\eps\in \Sc^\infty_\F$, $\int_0^.Z^\eps dW\in \textrm{BMO}(\P)$, we get from \textbf{(H$S$)} that $\int_0^.L^\eps dW\in \textrm{BMO}(\P)$.  Therefore, the process 
 $\Ec(\int_0^. L^\epsilon_s dW_s)$ is an $\F$-martingale from Theorem 2.3 in \cite{kaz94}. Applying the Girsanov theorem we get that the process $\bar W$ defined by 
\beqs
\bar W _t & := & W _t +\int_0^tL^\eps_sds\;,
\enqs 
for all $t\in[0,T]$, is a Brownian motion under the probability $\Q$ defined by 
\beqs
\frac{d \Q}{d \P} \Big|_{\Fc_T}  & = & \Ec \Big(-\int_0^T L^\epsilon_s dW_s \Big)\;.
\enqs
We also notice that under $\Q$, $(Y^\eps, Z^\eps)$ is solution to
\beq\label{EDSR eps ss Q}
Y^\eps_t & = & 1+\int_t^T \Big\{\frac{\lambda_s |\beta_s|^2}{|\sigma_s|^4} |\mu_s - \lambda_s \beta_s|^2  - \frac{|\mu_s - \lambda_s \beta_s|^2}{|\sigma_s|^2} Y^\eps_s  - 2 \lambda_s \beta_s \frac{(\mu_s - \lambda_s \beta_s)}{|\sigma_s|^2}  +\lambda_s \nonumber\\
 & &  - \;\lambda_s Y^\eps_s
- \frac{\big|\lambda_s \beta_s + (\lambda_s \beta_s - \mu_s) \frac{\lambda_s |\beta_s|^2}{|\sigma_s|^2}\big|^2}{|\sigma_s|^2 (Y^\eps_s \vee \eps) + \lambda_s |\beta_s|^2} \Big\}ds - \int_t^T Z^\eps_s d \bar W_s\;,\quad t\in[0,T]\;.
\enq

\vspace{2mm}

\ni\textbf{Step 4}: \textit{Comparison under the new probability measure $\Q$.}

 We first notice that the generator $\bar f^\epsilon$ of the BSDE \reff{EDSR eps ss Q} admits the following lower bound 
 \beqs 
 \bar f^\epsilon( t ,y,z)&  \geq & 
   \frac{\lambda_t |\beta_t|^2}{|\sigma_t|^4} |\mu_t - \lambda_t \beta_t|^2  + \lambda_t - \lambda_t y  - 2 \lambda_t \beta_t \frac{(\mu_t - \lambda_t \beta_t)}{|\sigma_t|^2}\\
  & &    -\frac{|\mu_t - \lambda_t \beta_t|^2}{|\sigma_t|^2} y
- \frac{\big|\lambda_t \beta_t + (\lambda_t \beta_t - \mu_t) \frac{\lambda_t |\beta_t|^2}{|\sigma_t|^2}\big|^2}{\lambda_t |\beta_t|^2} \1_{\lambda_t \beta_t \neq 0} \\ & = & 
 - \lambda_t y - \frac{|\mu_t - \lambda_t \beta_t|^2}{|\sigma_t|^2 }y \;,
\enqs
for all $(t,y,z)\in[0,T]\times\R\times\R$.

We now study the following BSDE
\beq\label{lowbdbarf}
\underline Y_t & = &  1+ \int_t^T \displaystyle{\Big[-\lambda_s - \frac{|\mu_s - \lambda_s \beta_s|^2}{|\sigma_s|^2 }\Big] \underline Y_s ds   - \int_t^T \underline Z_s d \bar W_s}\;,\quad t\in[0,T]\;.
\enq
Since this BSDE is linear, it has a unique solution given by (see e.g. \cite{nekpenque97})
\beqs
\underline Y_t & := & \E_\Q \Big[ \exp\Big( - \int_t^T \big( \lambda_s + \frac{|\mu_s - \lambda_s \beta_s|^2}{|\sigma_s|^2} \big)ds \Big) \Big| \Fc_t \Big]\;,\quad t\in[0,T]\;.
\enqs
Applying Theorem 2.2  of \cite{nekpenque97} for BSDEs \reff{EDSR eps ss Q} and \reff{lowbdbarf} we have
\beqs
Y^\epsilon_t & \geq & \underline Y_t\;,\quad t\in[0,T]\;.
\enqs 
By \reff{condeps}, we have $  \eps \leq \underline Y_t$ for any $t \in [0,T]$. Consequently, 
 $ Y^\eps_t \geq \eps$ for any $t \in [0,T]$, and $(Y^\eps,Z^\eps)$ is solution to \reff{EDSR0}. 
\ep

\vspace{2mm}

\ni We now are able to prove that the BSDE $(\mathfrak{f},1)$ admits a solution.
\begin{Proposition}\label{prop extce YZ}
The BSDE \reff{EDSRf^1} admits a solution $(Y,Z,U)\in \Sc^\infty_\G\times L^2_\G\times L^2(\lambda)$ with $Y\in \Sc^{\infty,+}_\G$. 
\end{Proposition}

\ni \textbf{Proof.}
From Theorem \ref{Thm-exist-gen} and Proposition \ref{sol f1avsaut}, we obtain that the BSDE \reff{EDSRf^1} admits a solution $(Y,Z,U)\in \Sc^\infty_\G\times L^2_\G\times L^2(\lambda)$, with $Y$ given by 
\beqs
Y_t & = & Y_t^b\mathds{1}_{\tau<t}+ \mathds{1}_{\tau\geq t}\;,\quad t\in[0,T]\;.
\enqs
 with $Y^b\in \Sc^{\infty,+}_\F$ from Proposition  \ref{sol f1avsaut}. Therefore $Y \in \Sc^{\infty,+}_\G$.
\ep
\vspace{2mm}

\subsection{Solution to the BSDE $(\mathfrak{g},H)$}

\vspace{2mm}

We first notice that the BSDE $(\mathfrak{g}, H)$ can be rewritten under the form
\begin{equation}\label{EDSR ronde}
\left\{\begin{array}{rcl}
d\Yc_t & = &  \displaystyle{\Big\{ \frac{(\mu_t Y_t + \sigma_t Z_t + \lambda^\G_t \beta_t U_t)( \sigma_t Y_t \Zc_t + \lambda^\G_t \beta_t (U_t + Y_t ) \Uc_t )}{Y_t (|\sigma_t|^2 Y_t + \lambda^\G_t |\beta_t|^2 (U_t + Y_t))} - \frac{Z_t}{Y_t} \Zc_t }\\ && ~~- \frac{\lambda^\G_t U_t}{Y_t} \Uc_t  -\; \lambda^\G_t \Uc_t \Big\} dt+  \Zc_t dW_t + \Uc_t dN_t\;, \quad t\in[0, T \wedge \tau] \;,\\
\Yc_{T \wedge \tau} & = & H \;.
\end{array}\right.
\end{equation}
Since $Y_t\mathds{1}_{t< \tau} = Y^b_t\mathds{1}_{t< \tau}$ and $U_t\mathds{1}_{t\leq \tau} = (1-Y^b_t)\mathds{1}_{t\leq \tau}$, we consider the associated decomposed BSDE in $\F$: find $(\Yc^b,\Zc^b)\in \Sc^\infty_\F\times L^2_\F$ such that 
\begin{equation}\label{EDSR0 ronde}
\left\{\begin{array}{rcl}
d\Yc^b_t & = &  \displaystyle{\Big\{ \frac{((\mu_t-\lambda_t \beta_t) Y^b_t + \sigma_t Z^b_t + \lambda_t \beta_t )( \sigma_t Y^b_t \Zc^b_t + \lambda_t \beta_t H^a_t- \lambda_t \beta_t  \Yc^b_t )}{Y^b_t (|\sigma_t|^2 Y^b_t + \lambda_t |\beta_t|^2 )} }\\ && ~ ~- \frac{Z^b_t}{Y^b_t} \Zc^b_t-\; \displaystyle{\frac{\lambda_t }{Y^b_t}  H^a_t+ \frac{\lambda_t }{Y^b_t} \Yc^b_t\Big\} dt+  \Zc^b_t dW_t} \;,\quad t \in[0,T]\;,\\
\Yc^b_{T} & = & H^b\;.
\end{array}\right.
\end{equation}
We notice that this BSDE has a Lipschitz generator w.r.t. the unknown $(\Yc^b,\Zc^b)$. However the Lipschitz coefficient depends on $Z^b$ which is not necessarily bounded. Thus we cannot apply the existing results and have to deal with this issue. 

\begin{Proposition}\label{sol f2avsaut}
The BSDE \reff{EDSR0 ronde} admits a solution $(\Yc^b, \Zc^b)$ in $\Sc^\infty_\F  \times L^2_\F$  with $\int_0^.\Zc^b dW\in \rm{BMO}(\P)$.
\end{Proposition}
\ni\textbf{Proof.}
We first define the equivalent probability $\Q$ to $\P$ defined by its Radon-Nikodym density $\frac{d\Q}{d\P} \big|_{\Fc_T} = \Ec (\int_0^T\rho_t dW_t)$ where $\rho$ is given by
\beqs
\rho_t &:=&  \frac{Z^b_t}{Y^b_t} - \frac{\sigma_t \big((\mu_t - \lambda_t \beta_t) Y^b_t + \sigma_t Z^b_t + \lambda_t \beta_t\big)}{|\sigma_t|^2 Y^b_t + \lambda_t |\beta_t|^2 }\;,\quad t\in[0,T]\;.
\enqs
Since $\int_0^.Z^b dW\in \mathrm{BMO}(\P)$, $Y^b \in \Sc^{\infty,+}_\F$ and the coefficients $\mu$, $\sigma$ and $\beta$ satisfy \textbf{(H$S$)}, it implies that $\int_0^.\rho dW \in \mathrm{BMO}(\P)$. 
Therefore, $\bar W_t := W_t - \int_0^t  \rho_s ds$ is a $\Q$-Brownian motion. Hence, the BSDE \reff{EDSR0 ronde} can be written as
\begin{equation}\label{equation YC lineaire}
\left\{\begin{array}{rcl}
d\Yc^b_t & = & a_t (\Yc^b_t - H^a_t) dt+ \Zc^b_t d\bar W_t \;,\quad t\in[0,T]\;,\\
\Yc^b_{T \wedge \tau} & = & H^b \;,
\end{array}\right.
\end{equation}
with
\beqs
a_t&:=&  \frac{\lambda_t |\sigma_t|^2 Y^b_t- \lambda_t \beta_t  ((\mu_t-\lambda_t \beta_t) Y^b_t + \sigma_t Z^b_t ) }{Y^b_t (|\sigma_t|^2 Y^b_t + \lambda_t |\beta_t|^2 )} \;,\quad t\in[0,T]\;.\\
\enqs
By definition of $a$ we can see that $\int_0^.a dW\in  \textrm{BMO}(\P)$ since the coefficients $\mu$, $\sigma$, $\beta$ and $\lambda$ are bounded, $Y^b \in \Sc^{\infty, +}_\F$ and $\int_0^.Z^b dW \in \rm{BMO}(\P)$.
Using Theorem \ref{Stabilite} with $\Q_1=\P$ and $\Q_2=\Q$, we get $\int_0^.a d\bar W \in \rm{BMO}(\Q)$. Therefore, there exists a constant $l' \geq 0$ such that $\E_\Q [ \int^T_\nu |a_s|^2  ds | \Fc_\nu] \leq l' $ for any $\nu \in \Tc_\F[0,T]$. We now prove that the process $\Yc^b$ defined by 
\beqs
\Yc^b_t & := & \E_\Q \Big[ \frac{\Gamma_T}{\Gamma_t}H^b + \int_t^T \frac{\Gamma_s}{\Gamma_t} a_s H^a_s ds \Big| \Fc_t \Big] \;,\quad t\in[0,T]\;,
\enqs
with $\Gamma_t := \exp(-\int_0^t a_s ds)$, is solution of the BSDE \reff{EDSR0 ronde}. We proceed in four steps. 

\vspace{2mm}

\ni\textbf{Step 1.} \textit{Integrability of the process $\Gamma$}.

We first prove that for any $p \geq 1$ there exists a constant $C > 0$ such that the process $\Gamma$ satisfies for any $t \in [0,T]$
\beq\label{cond int Gamma}
\E_\Q \Big[ \sup_{t \leq s \leq T} \Big| \frac{\Gamma_s}{\Gamma_t}\Big|^p \Big| \Fc_t \Big] & \leq & C\;.
\enq
Since $\E_\Q [ \int^T_\nu |a_s|^2  ds | \Fc_\nu] \leq l'$ for any $\nu\in \Tc_\F[0,T]$, we get from Proposition \ref{expbound} that there exists a constant $\delta$ such that $0 < \delta < \frac{1}{l'}$ and
\beqs
\E_\Q \Big[ \exp\Big( \delta \int_\nu^T |a_s|^2 ds \Big) \Big| \Fc_\nu \Big] & \leq & \frac{1}{1 - \delta l'}\;.
\enqs 
We get for any $0 \leq t \leq s \leq T$ 
\beqs
\Big| \frac{\Gamma_s}{\Gamma_t}\Big|^p & \leq & \exp\Big( \int_t^s \big(\delta |a_r|^2 + \frac{p^2}{4 \delta} \big) dr \Big)\\
&\leq & \exp\Big( \frac{p^2}{4 \delta} T \Big) \exp\Big( \delta \int_0^T |a_r|^2 dr \Big)\;.
\enqs
Consequently, we get
\beqs
\E_\Q \Big[ \sup_{t \leq s \leq T} \Big| \frac{\Gamma_s}{\Gamma_t}\Big|^p \Big| \Fc_t\Big] &\leq & \exp\Big( \frac{p^2}{4 \delta}T \Big)  \frac{1}{1 - \delta l'}\;.
\enqs

\vspace{2mm}

\ni\textbf{Step 2.} \textit{Uniform boundedness of $\Yc^b$}.

We now prove that $\Yc^b \in \Sc^\infty_\F$. For that we remark that by definition of $\Yc^b$ we have the following inequality
\beqs
|\Yc^b_t| &\leq & \|H^b\|_\infty \E_\Q \Big[ \frac{\Gamma_T}{\Gamma_t} \Big| \Fc_t \Big] + \|H^a\|_{\Sc^\infty} \E_\Q\Big[ \int_t^T |a_s|^2 ds \Big| \Fc_t \Big] + \|H^a\|_{\Sc^\infty}\E_\Q \Big[ \int_t^T   \Big|\frac{\Gamma_s}{\Gamma_t}\Big|^2 ds \Big| \Fc_t \Big]\;.
\enqs
Therefore, we get that $\Yc^b \in \Sc^\infty_\F$.

\vspace{2mm}

\ni \textbf{Step 3.} \textit{Dynamics of $\Yc^b$}.

 We now prove that $\Yc^b$ satisfies \reff{equation YC lineaire}. For that we introduce the $\Q$-martingale $m$ defined by 
 \beqs
 m_t & := & \Gamma_t \Yc^b_t + \int_0^t \Gamma_s a_s H^a_s ds\;,\quad t\in[0,T]\;.
 \enqs 
 We first notice that $m$ is $\Q$-square integrable. Indeed, from the definition of $m$, there exists a constant $C$ such that
 \beqs
 \E_\Q\Big[|m_t|^2\Big] & \leq & C\Big( \E_\Q\Big[\big| \Gamma_t \Yc^b_t \big|^2\Big]+\E_\Q\Big[\int_0^t\big| \Gamma_sa_sH^a_s \big|^2ds\Big] \Big)\;,
 \enqs
 for all $t\in[0,T]$. Since $\Yc^b \in \Sc^\infty_\F$, we get from \reff{Habded} 
 and from Cauchy-Schwarz inequality the existence of a constant $C$ such that
\beqs
 \E_\Q\Big[|m_t|^2\Big] & \leq & C\Big( \E_\Q\Big[\big| \Gamma_t  \big|^2\Big]+\sqrt{\E_\Q\Big[\Big(\int_0^t\big| a_s \big|^2ds\Big)^2\Big]}\sqrt{\E_\Q\Big[\sup_{0 \leq s \leq t}\big| \Gamma_s \big|^4\Big]} ~\Big)\;,
\enqs 
 for all $t\in[0,T]$. Since $\int_0^.a dW\in \rm{BMO} (\P)$ we have  from Theorem \ref{Stabilite}  $\int_0^.a d\bar W \in \rm{BMO}(\Q)$, and we get from Proposition \ref{expbound} and \reff{cond int Gamma}  
 \beqs
 \E_\Q\Big[|m_t|^2\Big] & < & \infty\;, \quad t\in[0,T]\;. 
 \enqs
Therefore, there exists a predictable process $\tilde \Zc$ such that $\E_\Q [\int_0^T | \tilde \Zc _s |^2 ds ] < \infty$ and
\beqs
\Gamma_t \Yc^b_t + \int_0^t \Gamma_s a_s H^a_s ds & = & m_0 + \int_0^t \tilde \Zc_s d \bar W_s\;,\quad t\in[0,T] \;.
\enqs
From It\^o's formula and the definition of $\Yc^b_T$ we have 
\beq\label{dynyb}
 \Yc^b_t & = &H^b-\int_t^T a_s (\Yc^b_s - H^a_s) ds - \int_t^T\Zc^b_t d \bar W_s \;,\quad t\in[0,T]\;.
\enq
where the process $\Zc^b$ is defined by 
\beqs
\Zc^b_t & := & \frac{\tilde \Zc_t}{\Gamma_t}\;,\quad t\in[0,T]\;.
\enqs
We now prove that $\int_0^.\Zc^b d\bar W\in\rm{BMO}(\Q)$. Using \reff{dynyb}, there exists a constant $C$ such that  
\beqs
\sup_{\nu\in \Tc_\F[0,T]}\E_\Q\Big[ \int_\nu^T|\Zc^b_s|^2ds \Big| \Fc_\nu \Big] & \leq & C\Big(( {\|\Yc^b\|}^2_{\Sc^\infty}+ {\|H^a\|}^2_{\Sc^\infty})\sup_{\nu\in \Tc_\F[0,T]}\E_\Q\Big[\int_\nu^T|a_s|^2ds \Big| \Fc_\nu \Big]   \\
& &+ {\|H^b\|}^2_{\infty} + {\|\Yc^b\|}^2_{\Sc^\infty}\Big)\;.
\enqs
Using $\Yc^b\in \Sc^\infty_\F$, \reff{Habded} and $\int_0^.a d\bar W\in \rm{BMO}(\Q)$, we get that $\int_0^.\Zc^b d\bar W \in \rm{BMO}(\Q)$. Thus, using ${d\P \over d\Q}\big|_{\Fc_T}=\Ec(-\int_0^.\rho d\bar W)_T$ and Theorem \ref{Stabilite} with $\Q_1=\Q$ and $\Q_2=\P$  we obtain that  
\beqs
\int_0^.\Zc^b d W~~=~~\int_0^.\Zc^b d\bar W -\langle \int_0^.\Zc^b d\bar W, \int_0^.\rho d\bar W\rangle & \in & \rm{BMO}(\P)\;.
\enqs
To conclude we get from \reff{dynyb} and the definition of $\bar W$ that $(\Yc^b, \Zc^b)$ is a solution to the BSDE \reff{EDSR0 ronde}.\ep

\vspace{2mm}

\ni We now prove the existence of a solution to the BSDE $(\mathfrak{g},H)$.

\begin{Proposition}\label{prop extce YcZc}
The BSDE \reff{EDSRf^2} admits a solution $(\Yc,\Zc,\Uc)\in \Sc^\infty_\G\times L^2_\G\times L^2(\lambda)$. 
\end{Proposition}

\ni \textbf{Proof.} From Theorem \ref{Thm-exist-gen} and Proposition \ref{sol f2avsaut}, we obtain that the BSDE \reff{EDSRf^2} admits a solution $(\Yc,\Zc,\Uc)\in \Sc^\infty_\G\times L^2_\G\times L^2(\lambda)$. \ep

\subsection{Solution to the BSDE $(\mathfrak{h},0)$} 
 We recall that the BSDE $(\mathfrak{h},0)$ is
 \beq \label{EDSR 3}
 \Upsilon_t & = & \int_{t\wedge \tau}^{T \wedge \tau} \Big(|\Zc_s|^2 Y_s + \lambda^\G_s (U_s + Y_s) |\Uc_s|^2
\displaystyle{-\frac{|\sigma_sY_s\Zc_s+\lambda^\G_s\beta_s\Uc_s(U_s+Y_s)|^2}{|\sigma_s|^2 Y_s + \lambda^\G_s |\beta_s|^2 ( U_s + Y_s)}} \Big) ds\nonumber \\
 & &  - \int_{t\wedge \tau}^{T \wedge \tau} \Xi_s dW_s - \int_{t\wedge \tau}^{T \wedge \tau} \Theta_s dM_s\;,\quad t\in[0,T]\;.
 \enq
Using  the definitions of $Y$, $U$, $\Zc$ and $\Uc$, we therefore consider the associated decomposed BSDE in $\F$: find $(\Upsilon^b,\Xi^b) \in \Sc^\infty_\F \times L^2_\F$ such that
 \beq
 \Upsilon_t^b & = & \int_{t}^{T} \Big(|\Zc_s^b|^2 Y^b_s + \lambda_s | H^a_s-\Yc^b_s|^2
\displaystyle{-\frac{|\sigma_sY^b_s\Zc^b_s+\lambda_s\beta_s (H^a_s-\Yc^b_s)|^2}{|\sigma_s|^2 Y^b_s + \lambda_s |\beta_s|^2}} -\lambda_s \Upsilon_s^b\Big) ds\nonumber \\\label{EDSR3decom}
 & & 
  - \int_{t}^{T} \Xi_s^b dW_s \;,\qquad t\in[0,T]\;.
  \enq
   \begin{Proposition}\label{sol f3avsaut}
 The BSDE \reff{EDSR3decom} admits a solution $(\Upsilon^b , \Xi^b)\in \Sc^\infty_\F\times L^2_\F$. 
 \end{Proposition} 
\ni \textbf{Proof.} Denote by $R$ the process defined by
 \beqs
 R_t & := & |\Zc_t^b|^2 Y^b_t + \lambda_t  | H^a_t-\Yc^b_t|^2
\displaystyle{-\frac{|\sigma_tY^b_t\Zc^b_t+\lambda_t\beta_t (H^a_t-\Yc^b_t)|^2}{|\sigma_t|^2 Y^b_t + \lambda_t |\beta_t|^2 }}  \;,
 \enqs
 for $t\in[0,T]$.
Define the process $\tilde \Upsilon^b$ by 
\beqs
\tilde \Upsilon^b_t & := & \E\Big[ \int_t^TR_se^{-\int_0^s\lambda_udu}ds \Big|\Fc_t\Big]\;,\quad t\in[0,T]\;.
\enqs
From \textbf{(H$S$)}, $\lambda$ is bounded, $Y^b \in \Sc^{\infty,+}_\F$, $H^a \in \Sc^\infty_\F$, $\Yc^b \in \Sc^\infty_\F$ and $\int_0^.\Zc^b dW\in \mathrm{BMO}(\P)$, we get from Proposition \ref{expbound} that
$\tilde \Upsilon^b\in \Sc^\infty_\F$ and the process $\tilde \Upsilon^b+\int_0^.R_se^{-\int_0^s\lambda_udu}ds$ is a square integrable martingale. Hence there exists a process $\tilde \Xi^b\in L^2_\F$ such that 
 \beqs
 \tilde \Upsilon^b_t & = & \int_t^TR_se^{-\int_0^s\lambda_udu}ds -\int_t^T\tilde \Xi^b_sdW_s\;,\quad t\in[0,T]\;. 
 \enqs
  From It\^o's formula we get that the processes $(\Upsilon^b,\Xi^b)$ defined by 
\beqs
\Upsilon^b_t ~=~ \tilde \Upsilon^b_t e^{\int_0^t\lambda_sds}   & \mbox{ and } & \Xi^b_t ~=~ \tilde \Xi^b_t e^{\int_0^t\lambda_sds} 
\enqs
satisfy \reff{EDSR3decom}. Since $\tilde \Xi^b\in L^2_\F$ and $\lambda$ is uniformly bounded we get that $\Xi^b\in L^2_\F$. Finally, since $\tilde \Upsilon^b\in \Sc^\infty_\F$ we get that $\Upsilon^b\in \Sc^\infty_\F$. 
 \ep

\vspace{2mm}

Finally, we prove the existence of a solutin to the BSDE $(\mathfrak{h},0)$.
\begin{Proposition}
The BSDE \reff{EDSRf^3} admits a solution $(\Upsilon,\Xi,\Theta)\in \Sc^\infty_\G\times L^2_\G\times L^2(\lambda)$. 
\end{Proposition}

\ni \textbf{Proof.} From Theorem \ref{Thm-exist-gen} and Proposition \ref{sol f3avsaut}, we obtain that the BSDE \reff{EDSRf^3} admits a solution $(Y,Z,U)\in \Sc^\infty_\G\times L^2_\G\times L^2(\lambda)$.
\ep
\appendix

\section{Appendix}
\setcounter{equation}{0} \setcounter{Assumption}{0}
\setcounter{Theorem}{0} \setcounter{Proposition}{0}
\setcounter{Corollary}{0} \setcounter{Lemma}{0}
\setcounter{Definition}{0} \setcounter{Remark}{0}

\subsection{Proof of Proposition \ref{propjeulin}}
We first suppose that $X$ is a nonnegative $\Pc(\G)$-measurable process.  For $n\geq 1$, we define the process $X^n$ by
\beqs
X^n_t & = & X_t\wedge n\;,\quad t\in [0,T]\;.
\enqs 
Then $X^n$ is a bounded $\G$-predictable process, and from Lemma 4.4 in \cite{jeu80}, there exist a $\Pc(\F)$-measurable process  $X^{n,b}$ and a  $\Pc(\F)\otimes\Bc(\R_+)$-measurable process $X^{n,a}$ such that
\beq\label{decompXn}
X_t^n & = & X^{n,b}_t \1_{t \leq \tau }+ X^{n,a}_t(\tau) \1_{t > \tau} \;,\quad t\in [0,T]\;.
\enq
Since the sequence $(X^n)_n$ is nondecreasing, we can assume w.l.o.g. that the sequences $(X^{a,n})_n$ and $(X^{b,n})_n$ are also nondecreasing. Define the processes $X^a$ and $X^b$ by
\beqs
X^a ~ = ~ \lim_{n\rightarrow\infty} X^{n,a} & \mbox{ and } & X^b ~ = ~ \lim_{n\rightarrow\infty} X^{n,b} \;.
\enqs
Then $X^a$ is $\Pc(\F)\otimes\Bc(\R_+)$-measurable and $X^b$ is $\Pc(\F)$-measurable and sending $n$ to infinity in \reff{decompXn}, we get
\beq\label{decompX}
X_t & = & X^{b}_t \1_{t \leq \tau }+ X^{a}_t(\tau) \1_{t > \tau} \;,\quad t\in [0,T]\;.
\enq
For a general $\Pc(\G)$-measurable process $X$, we write $X=X^+-X^-$ where $X^+=\max(X,0)$ and $X^-=\max(-X,0)$ and we apply the previous result to the nonnegative processes $X^+$ and $X^-$. From the linear stability of the decomposition \reff{decompX} we get the result. \ep

\subsection{BMO Stability}

\begin{Theorem} \label{Stabilite} Let $\Q_1$ and $\Q_2$ be two probability measures on $(\Omega,\Gc)$. Let $M$ and $N$ be two continuous $(\F,\Q_1)$-local martingales with $N\in\rm{BMO}(\mathbb{Q}_1)$.  Suppose that $\Q_1$ and $\Q_2$ are equivalent with ${d\mathbb{Q}_2\over d\mathbb{Q}_1}\big|_{\Fc_T}={\mathcal E}(N)_T$. If $M\in\rm{BMO}(\mathbb{Q}_1)$ then $M-\langle M,N\rangle\in\rm{BMO}(\mathbb{Q}_2)$.  
\end{Theorem}
\ni\textbf{Proof.} This result is a direct consequence of  Theorem 3.6 in   \cite{kaz94}.
\ep
\subsection{An estimate for conditional moments}
\begin{Proposition} \label{expbound} Let $A$ be a continuous increasing $\F$-adapted process.  Fix a $t\geq 0$ such that there exists a constant $C > 0$ satisfying 
\beqs
\mathbb{E}\big[A_t-A_s \big| \mathcal{F}_s\big] & \leq & C \;,
\enqs
for  any $ s \in [0, t]$.
Then, we have for any $s \in [0,t]$ and any $p \geq 1$
\beqs
 \mathbb{E} \big[{|A_t-A_s|}^{p} | {\mathcal F}_s \big]\leq p! |C|^{p}
\enqs
and
\beqs
\mathbb{E}\Big[\exp \big(\delta(A_t-A_s) \big) \big| \mathcal{F}_s\Big] & \leq & \frac{1}{1-\delta C} \;,
\enqs
for any $\delta\in (0,{1\over C})$.
\end{Proposition}

\ni\textbf{Proof.}   Let $A$ be a continuous increasing $\F$-adapted process satisfying $\mathbb{E}[A_t-A_s | \mathcal{F}_s] \leq C$ for any $s \in[0,t ]$. We first prove by iteration that $\mathbb{E} [{|A_t-A_s|}^{p} | {\mathcal F}_s]\leq p! |C|^{p}$ for any $p \geq 1$. 

\vspace{2mm}

\ni $\bullet$ For $p=1$, we have by assumption
$\mathbb{E}[A_t-A_s | \mathcal{F}_s] \leq C$.

\vspace{2mm}

\ni $\bullet$ 
Suppose that for some $p\geq 2$, we have $\mathbb{E}[{|A_t-A_s|}^{p-1} | {\mathcal F}_s]\leq (p-1)! |C|^{p-1}$. Since $A$ is a continuous increasing $\F$-adapted process we have 
\beqs
{|A_t-A_s|}^{p}=p \int_s^t {|A_t-A_u|}^{p-1} dA_u\;,\enqs
 for any $s\in[0, t]$. Consequently we get
\beqs
\E \big[|A_t-A_s|^{p} \big| \Fc_s \big] & = & p \E \Big[ \int_s^t |A_t-A_u|^{p-1} dA_u \Big| \Fc_s \Big] \\
& = & p \E \Big[ \int_s^t \E \Big[|A_t-A_u|^{p-1} \Big| \Fc_u \Big] dA_u \Big| \Fc_s \Big] \\
& \leq & p! |C|^{p-1} \E [ A_t - A_s | \Fc_s] \\
& \leq & p! |C|^p \;.
\enqs
\ni$\bullet$ Since the result holds true for $p=1$ and for any $p\geq 2$ as soon as it holds for $p-1$, it holds for $p$, we get 
\beqs
\E\big[{|A_t-A_s\big|}^{p} \big| {\Fc}_s]\leq p! |C|^{p}\;,
\enqs
for any $p\geq 1$.

\vspace{2mm}

\ni From this last inequality, we get for any $\delta\in (0, \frac{1}{C})$
\beqs
\E\Big[\sum_{p\geq 0}{1\over p!}|\delta|^p {|A_t-A_s|}^p\Big| {\mathcal F_s}\Big] & \leq &  \sum_{p\geq 0} |\delta C|^p~=~{1\over 1-\delta C}\;,
\enqs
which is the expected result. \ep

%

\end{document}